\documentclass [a4paper]{article}
\usepackage[cp1251]{inputenc}
\usepackage{amsmath,amssymb}
\usepackage{cite}
\begin{document}

\title{Determinantal representations of   W-weighted Drazin inverse solutions of some quaternion matrix equations}
\author {{\textbf{Ivan I. Kyrchei}} \thanks{E-mail address: kyrchei@online.ua}
\\Pidstrygach Institute for Applied Problems of Mechanics and\\ Mathematics NAS of Ukraine, Lviv, 79060, Ukraine}

\date{}
 \maketitle

\textbf{Abstract:}
By using  determinantal representations of  the W-weighted Drazin inverse  previously obtained by the author within the framework of the theory of the column-row determinants, we get explicit formulas for determinantal representations of the W-weighted Drazin inverse solutions
 (analogs of Cramer's rule) of the quaternion matrix equations $ {\bf W}{\bf A}{\bf W}{\bf X}={\bf D}$,  $ {\bf X}{\bf W}{\bf A}{\bf W}={\bf D} $, and ${\bf W}_{1}{\bf A}{\bf W}_{1}{\bf X}{\bf W}_{2}{\bf B}{\bf W}_{2}={\bf D} $.



\newtheorem{corollary}{Corollary}[section]
\newtheorem{theorem}{Theorem}[section]
\newtheorem{lemma}{Lemma}[section]
\newtheorem{definition}{Definition}[section]
\newtheorem{remark}{Remark}[section]
\newcommand{\rank}{\mathop{\rm rank}\nolimits}
\newtheorem{proposition}{Proposition}[section]

\section{Introduction}

Throughout the paper, we denote the real number field by ${\rm
{\mathbb{R}}}$, the set of all $m\times n$ matrices over the
quaternion algebra
\[{\rm {\mathbb{H}}}=\{a_{0}+a_{1}i+a_{2}j+a_{3}k\,
|\,i^{2}=j^{2}=k^{2}=-1,\, a_{0}, a_{1}, a_{2}, a_{3}\in{\rm
{\mathbb{R}}}\}\]
by ${\rm {\mathbb{H}}}^{m\times n}$, and by ${\rm {\mathbb{H}}}^{m\times n}_{r}$ the set of all $m\times n$ matrices over $\mathbb{H}$
with a rank $r$.  Let ${\rm M}\left( {n,{\rm {\mathbb{H}}}} \right)$ be the
ring of $n\times n$ quaternion matrices.
 For ${\rm {\bf A}}
 \in {\rm {\mathbb{H}}}^{n\times m}$, the symbols ${\rm {\bf A}}^{ *}$ stands for the conjugate transpose (Hermitian adjoint) matrix
of ${\rm {\bf A}}$.
 The matrix ${\rm {\bf A}} = \left( {a_{ij}}  \right) \in {\rm
{\mathbb{H}}}^{n\times n}$ is Hermitian if ${\rm {\bf
A}}^{ *}  = {\rm {\bf A}}$.

In the past, researches into the quaternion skew field  had more a theoretical importance, but now a growing number of investigations   give wide practical applications of quaternions.
In particular through their attitude orientation, the quaternions arise in various fields such as quaternionic quantum theory \cite{ad}, fluid mechanics and particle dynamics \cite{gi, gi1}, computer graphics \cite{ha}, aircraft orientation \cite{st}, robotic systems \cite{pe},       life science \cite{is,pro} and etc.

Research on quaternion matrix equations and  generalized inverses, which are usefulness tools used to solve  matrix equations, has been actively ongoing for more recent years. We mention  only some recent papers.
 Yuan,  Wang and Duan \cite{wa1} derived  solutions of the quaternion matrix equation $AX=B$  and their applications in color image restoration. Wang and Yu \cite{wa2} studied extreme ranks of real matrices in solution of the
quaternion matrix equation $AXB = C$. Yuan,  Liao and Lei \cite{yu} obtained the expressions of least squares Hermitian solution with minimum norm of the quaternion matrix equation $(AXB,CXD)=(E,F)$.
Feng and Cheng \cite{fe} gave a clear description of the solution set to the
quaternion matrix equation $AX -  \bar{X}B = 0$.  Jiang and
Wei \cite{ji} derived the explicit solution of the quaternion matrix equation
$X -  A \tilde{X}B = C$. Song, Chen and Wang \cite{son} obtained the expressions of
the explicit solutions of quaternion matrix equations $XF -  AX = BY$ and $XF  - A \tilde{X} = BY$.
Yuan and Wang \cite{yu1} gave the expressions of the least squares
$\eta$-Hermitian solution with the least norm of the quaternion matrix equation
$AXB+CXD=E$.
 Zhang,  Wei,  Lia and Zhao derived \cite{zha} the expressions of the minimal norm least squares solution, the pure imaginary least squares solution, and the real least squares solution for the quaternion matrix equation $AX=B$.

The  definitions of the generalized inverse matrices have been  extended  to quaternion matrices as follows.

The Moore-Penrose inverse of ${\bf A}\in{\rm {\mathbb{H}}}^{m\times n}$, denoted by ${\bf A}^{\dagger}$, is the unique matrix ${\bf
X}\in{\rm {\mathbb{H}}}^{n\times m}$ satisfying the following equations
 \begin{gather*} 1)\, {\rm {\bf A}}{\bf X}
{\rm {\bf A}} = {\rm {\bf A}};  2)\,
                                  {\bf X} {\rm {\bf
A}}{\bf X}  = {\bf X}; 3)\,
                                  \left( {\rm {\bf A}}{\bf X} \right)^{ *}  = {\rm
{\bf A}}{\bf X};  4)\,
                                \left( {{\bf X} {\rm {\bf A}}} \right)^{ *}  ={\bf X} {\rm {\bf A}}. \end{gather*}
 For ${\bf A}\in{\rm {\mathbb{H}}}^{n\times n}$ with $k = Ind\,{\bf A}$ the smallest positive number such that $\rank {\bf
 A}^{k+1}=\rank {\bf A}^{k}$, the
Drazin inverse of ${\bf A}$, denoted by ${\bf A}^{D}$, is defined to be the unique matrix ${\bf X}$ that satisfying Eq. 2)
and the  equations
                                \begin{gather*}
                                  5)\, {\rm {\bf A}}{\bf X}  ={\bf X} {\rm
{\bf A}}; 6)\,
                                   {\rm {\bf A}}^{k+1}{\rm
{\bf X}}={\rm {\bf
         A}}^{k}. \end{gather*}
 In particular, when $Ind{\kern 1pt} {\rm {\bf A}}=1$,
then ${\rm {\bf X}}$  is called the group inverse of ${\bf A}$ and is denoted by ${\rm {\bf X}}={\rm {\bf A}}^{g }$.
If $Ind{\kern 1pt} {\rm {\bf A}}=0$, then ${\rm {\bf A}}$ is
nonsingular, and ${\rm {\bf A}}^{D}\equiv {\bf A}^{\dagger}= {\rm {\bf A}}^{-1}$.

Cline and Greville \cite{cl} extended the Drazin inverse of square matrix to rectangular matrix,
 that has been generalized to the quaternion algebra as follows.

For ${\bf A}\in{\rm {\mathbb{H}}}^{m\times n}$ and ${\bf W}\in{\rm {\mathbb{H}}}^{n\times m}$, the W-weighted Drazin inverse of ${\bf
A}$ with respect to ${\bf W}$ is  the unique solution to the following equations
\begin{gather*}
7)\,  ({{\bf A}{\bf W}})^{k+1}{\rm
{\bf X}}{\bf W}=({\rm {\bf
         A}}{\bf W})^{k}; 8)\,
  {\rm {\bf X}}{\bf W}{\rm {\bf A}}{\bf W}{\rm {\bf X}}={\rm {\bf X}}; 9)\,
   {\rm
{\bf A}}{\bf W}{\rm {\bf X}}={\rm {\bf X}}{\bf W}{\rm {\bf A}},
\end{gather*}
where $k= {\rm max}\{Ind({\bf A}{\bf W}), Ind({\bf W}{\bf A})\}$.

The Drazin inverse and weighted
Drazin inverse has several important applications such as, applications in singular differential and difference equations\cite{ca}, signal processing \cite{camp}, Marckov chains and statistic problems\cite{ca1,sp},
descriptor continuous-time
systems \cite{ka},   numerical analysis and Kronecker
product systems \cite{ze}, solving singular fuzzy linear system \cite{nik,nik1}, constrained linear systems \cite{wang} and etc.

 Cramer's rule for the W-weighted Drazin inverse
solutions, in particular, have been derived in \cite{wei} for  singular linear equations and   in \cite{wang} for a class of restricted matrix equations.
Recently, within the framework of the theory of the column-row determinants  Song \cite{song_ela}  has first obtained a determinantal
representation of the
W-weighted Drazin inverse  and Cramer's rule of a class of restricted matrix equations over the quaternion algebra.
But in obtaining, he has used  auxiliary matrices other than that are given.
In  \cite{ky11}, we have obtained new determinantal representations of the W-weighted Drazin inverse over the quaternion skew field without any auxiliary matrices.

  An important application of  determinantal representations of  generalized  inverses is the  Cramer rule for   generalized
  inverse solutions of matrix equations.

But, when there is a need for a W-weighted Drazin inverse solution?
Consider for example the following matrix equation,
$
  {\rm {\bf A}}_{1}{\rm {\bf X}}={\rm {\bf                                                                                                                   D}}$.
If ${\bf A}_{1}$ is rectangular and we can represent it as
${\rm {\bf A}}_{1}= {\bf W}{\bf A}{\bf W}$, where ${\bf W}{\bf A}$ and ${\bf A}{\bf W}$ are quadratic and singular, then  its W-weighted Drazin inverse solution is needed.

     In the paper we  investigate analogs of  Cramer's rule for   W-weighted Drazin inverse solutions of the following matrix equations over the quaternion skew field ${\rm {\mathbb{H}}}$,
\begin{gather}\label{eq:1ax}
{\bf W}{\bf A}{\bf W}{\bf X}={\bf D},\\
\label{eq:1xa}
{\bf X}{\bf W}{\bf A}{\bf W}={\bf D},\\
\label{eq:1axb}
{\bf W}_{1}{\bf A}{\bf W}_{1}{\bf X}{\bf W}_{2}{\bf B}{\bf W}_{2}={\bf D}.
\end{gather}
The paper is organized as follows. We start with some
basic concepts and results from the theory of  the row-column determinants   in Subsection 2.1. In Subsection 2.2, we  give  determinantal representations of the Moore-Penrose and Drazin inverses for a  quaternion matrix.  Determinantal representations of the W-weighted Drazin inverse and its properties we consider in  Subsection 2.3.  In Subsection 3.1, we  give the background of the problem of  Cramer's rule for the W-weighted Drazin inverse solution. In Subsection 3.2 we
obtain explicit representation formulas of the W-weighted Drazin inverse solutions (analogs of Cramer's rule) of the quaternion matrix equations    (\ref{eq:1ax}-\ref{eq:1axb}). In Section 4, we give
 numerical examples to illustrate the main result.

\section{Preliminaries}
\subsection{Elements of the theory of the column and row determinants}
The theory of the row-column determinants over the quaternion skew field has been introduced in \cite{ky1,ky2, ky3}, and later it has been applied to research  generalized inverses and generalized inverse solutions of matrix equations. In particular, determinantal representations of the Moore-Penrose \cite{ky4,ky7} and explicit representation formulas for the minimum norm least squares solutions of some quaternion matrix equations \cite{ky6}, and determinantal representations of the   Drazin \cite{ky5}, and W-weighted Drazin inverses \cite{ky11} have been obtained by the  author. Song at al. derived determinantal representation of the generalized inverse $A_{T,S}^{2}$ \cite{song3}, Bott-Duffin inverse \cite{song5} and the Cramer rule for the  solutions of restricted matrix equations \cite{song1}, and the generalized Stein quaternion matrix equation   \cite{song6}, etc.

 For  ${\rm {\bf A}}=(a_{ij}) \in {\rm
M}\left( {n,{\mathbb{H}}} \right)$ we define $n$ row determinants and $n$ column determinants as follows.

Suppose $S_{n}$ is the symmetric group on the set $I_{n}=\{1,\ldots,n\}$.
\begin{definition}
 The $i$-th row determinant of ${\rm {\bf A}}=(a_{ij}) \in {\rm
M}\left( {n,{\mathbb{H}}} \right)$ is defined  for all $i = \overline{1,n} $
by putting
 \begin{gather*}{\rm{rdet}}_{ i} {\rm {\bf A}} =
{\sum\limits_{\sigma \in S_{n}} {\left( { - 1} \right)^{n - r}{a_{i{\kern
1pt} i_{k_{1}}} } {a_{i_{k_{1}}   i_{k_{1} + 1}}} \ldots } } {a_{i_{k_{1}
+ l_{1}}
 i}}  \ldots  {a_{i_{k_{r}}  i_{k_{r} + 1}}}
\ldots  {a_{i_{k_{r} + l_{r}}  i_{k_{r}} }},\\
\sigma = \left(
{i\,i_{k_{1}}  i_{k_{1} + 1} \ldots i_{k_{1} + l_{1}} } \right)\left(
{i_{k_{2}}  i_{k_{2} + 1} \ldots i_{k_{2} + l_{2}} } \right)\ldots \left(
{i_{k_{r}}  i_{k_{r} + 1} \ldots i_{k_{r} + l_{r}} } \right),\end{gather*}
with
conditions $i_{k_{2}} < i_{k_{3}}  < \ldots < i_{k_{r}}$ and $i_{k_{t}}  <
i_{k_{t} + s} $ for $t = \overline{2,r} $ and $s =\overline{1,l_{t}} $.
\end{definition}
\begin{definition}
The $j$-th column determinant
 of ${\rm {\bf
A}}=(a_{ij}) \in {\rm M}\left( {n,{\mathbb{H}}} \right)$ is defined for
all $j =\overline{1,n} $ by putting
 \begin{gather*}{\rm{cdet}} _{{j}}\, {\rm {\bf A}} =
{{\sum\limits_{\tau \in S_{n}} {\left( { - 1} \right)^{n - r}a_{j_{k_{r}}
j_{k_{r} + l_{r}} } \ldots a_{j_{k_{r} + 1} i_{k_{r}} }  \ldots } }a_{j\,
j_{k_{1} + l_{1}} }  \ldots  a_{ j_{k_{1} + 1} j_{k_{1}} }a_{j_{k_{1}}
j}},\\
\tau =
\left( {j_{k_{r} + l_{r}}  \ldots j_{k_{r} + 1} j_{k_{r}} } \right)\ldots
\left( {j_{k_{2} + l_{2}}  \ldots j_{k_{2} + 1} j_{k_{2}} } \right){\kern
1pt} \left( {j_{k_{1} + l_{1}}  \ldots j_{k_{1} + 1} j_{k_{1} } j}
\right), \end{gather*}
\noindent with conditions, $j_{k_{2}}  < j_{k_{3}}  < \ldots <
j_{k_{r}} $ and $j_{k_{t}}  < j_{k_{t} + s} $ for  $t = \overline{2,r} $
and $s = \overline{1,l_{t}}  $.
\end{definition}
Suppose ${\rm {\bf A}}_{}^{i{\kern 1pt} j} $ denotes the submatrix of
${\rm {\bf A}}$ obtained by deleting both the $i$th row and the $j$th
column. Let ${\rm {\bf a}}_{.j} $ be the $j$-th column and ${\rm {\bf
a}}_{i.} $ be the $i$-th row of ${\rm {\bf A}}$. Suppose ${\rm {\bf
A}}_{.j} \left( {{\rm {\bf b}}} \right)$ denotes the matrix obtained from
${\rm {\bf A}}$ by replacing its $j$-th column with the column ${\rm {\bf
b}}$, and ${\rm {\bf A}}_{i.} \left( {{\rm {\bf b}}} \right)$ denotes the
matrix obtained from ${\rm {\bf A}}$ by replacing its $i$-th row with the
row ${\rm {\bf b}}$.

The following theorem has a key value in the theory of the column and row
determinants.
\begin{theorem} \cite{ky1}\label{theorem:
determinant of hermitian matrix} If ${\rm {\bf A}} = \left( {a_{ij}}
\right) \in {\rm M}\left( {n,{\rm {\mathbb{H}}}} \right)$ is Hermitian,
then ${\rm{rdet}} _{1} {\rm {\bf A}} = \cdots = {\rm{rdet}} _{n} {\rm {\bf
A}} = {\rm{cdet}} _{1} {\rm {\bf A}} = \cdots = {\rm{cdet}} _{n} {\rm {\bf
A}} \in {\rm {\mathbb{R}}}.$
\end{theorem}
Since all  column and  row determinants of a
Hermitian matrix over ${\rm {\mathbb{H}}}$ are equal, we can define the
determinant of a  Hermitian matrix ${\rm {\bf A}}\in {\rm M}\left( {n,{\rm
{\mathbb{H}}}} \right)$. By definition, we put,
$\det {\rm {\bf A}}: = {\rm{rdet}}_{{i}}\,
{\rm {\bf A}} = {\rm{cdet}} _{{i}}\, {\rm {\bf A}}, $
 for all $i =\overline{1,n}$.
The determinant of a Hermitian matrix has properties similar to a usual determinant. They are completely explored
in
 \cite{ky1, ky2}
 by its row and
column determinants. They can be summarized by the following
theorems.
\begin{theorem}\label{theorem:row_combin} If the $i$-th row of
a Hermitian matrix ${\rm {\bf A}}\in {\rm M}\left( {n,{\rm
{\mathbb{H}}}} \right)$ is replaced with a left linear combination
of its other rows, i.e. ${\rm {\bf a}}_{i.} = c_{1} {\rm {\bf
a}}_{i_{1} .} + \ldots + c_{k}  {\rm {\bf a}}_{i_{k} .}$, where $
c_{l} \in {{\rm {\mathbb{H}}}}$ for all $ l = \overline{1, k}$ and
$\{i,i_{l}\}\subset I_{n} $, then
\[
 {\rm{rdet}}_{i}\, {\rm {\bf A}}_{i \, .} \left(
{c_{1} {\rm {\bf a}}_{i_{1} .} + \ldots + c_{k} {\rm {\bf
a}}_{i_{k} .}}  \right) = {\rm{cdet}} _{i}\, {\rm {\bf A}}_{i\, .}
\left( {c_{1}
 {\rm {\bf a}}_{i_{1} .} + \ldots + c_{k} {\rm {\bf
a}}_{i_{k} .}}  \right) = 0.
\]
\end{theorem}
\begin{theorem}\cite{ky1}\label{theorem:colum_combin} If the $j$-th column of
 a Hermitian matrix ${\rm {\bf A}}\in
{\rm M}\left( {n,{\rm {\mathbb{H}}}} \right)$   is replaced with a
right linear combination of its other columns, i.e. ${\rm {\bf
a}}_{.j} = {\rm {\bf a}}_{.j_{1}}   c_{1} + \ldots + {\rm {\bf
a}}_{.j_{k}} c_{k} $, where $c_{l} \in{{\rm {\mathbb{H}}}}$ for
all $ l = \overline{1, k}$ and $\{j,j_{l}\}\subset J_{n}$, then
 \[{\rm{cdet}} _{j}\, {\rm {\bf A}}_{.j}
\left( {{\rm {\bf a}}_{.j_{1}} c_{1} + \ldots + {\rm {\bf
a}}_{.j_{k}}c_{k}} \right) ={\rm{rdet}} _{j} \,{\rm {\bf A}}_{.j}
\left( {{\rm {\bf a}}_{.j_{1}}  c_{1} + \ldots + {\rm {\bf
a}}_{.j_{k}}  c_{k}} \right) = 0.
\]
\end{theorem}
  The determinant of a Hermitian matrix also has a property of expansion along arbitrary rows and columns using  row and column determinants of submatrices.

  We have the following theorem on the determinantal representation
of the inverse matrix over ${\mathbb{H}}$.
\begin{theorem}\cite{ky1} \label{theorem:deter_inver} The necessary and sufficient condition of invertibility
of  ${\rm {\bf A}} \in {\rm M}(n,{{\rm {\mathbb{H}}}})$ is
${\rm{ddet}} {\rm {\bf A}} \ne 0$. Then there exists ${\rm {\bf
A}}^{ - 1} = \left( {L{\rm {\bf A}}} \right)^{ - 1} = \left(
{R{\rm {\bf A}}} \right)^{ - 1}$, where
\begin{equation}
\label{eq:det_inv_rdet} \left( {L{\rm {\bf A}}} \right)^{ - 1}
=\left( {{\rm {\bf A}}^{ *}{\rm {\bf A}} } \right)^{ - 1}{\rm {\bf
A}}^{ *} ={\frac{{1}}{{{\rm{ddet}}{ \rm{\bf A}} }}}
\begin{pmatrix}
  {\mathbb{L}} _{11} & {\mathbb{L}} _{21}& \ldots & {\mathbb{L}} _{n1} \\
  {\mathbb{L}} _{12} & {\mathbb{L}} _{22} & \ldots & {\mathbb{L}} _{n2} \\
  \ldots & \ldots & \ldots & \ldots \\
 {\mathbb{L}} _{1n} & {\mathbb{L}} _{2n} & \ldots & {\mathbb{L}} _{nn}
\end{pmatrix},\end{equation}
\begin{equation}
\label{eq:det_inv_cdet}
 \left( {R{\rm {\bf A}}} \right)^{ - 1} = {\rm {\bf
A}}^{ *} \left( {{\rm {\bf A}}{\rm {\bf A}}^{ *} } \right)^{ - 1}
= {\frac{{1}}{{{\rm{ddet}}{ \rm{\bf A}}^{ *} }}}
\begin{pmatrix}
 {\mathbb{R}} _{11} & {\mathbb{R}} _{21} &\ldots & {\mathbb{R}} _{n1} \\
 {\mathbb{R}} _{12} & {\mathbb{R}} _{22} &\ldots & {\mathbb{R}} _{n2}  \\
 \ldots  & \ldots & \ldots & \ldots \\
 {\mathbb{R}} _{1n} & {\mathbb{R}} _{2n} &\ldots & {\mathbb{R}} _{nn}
\end{pmatrix}
\end{equation}
and \[{\mathbb{L}} _{ij} = {\rm{cdet}} _{j} ({\rm {\bf
A}}^{\ast}{\rm {\bf A}})_{.j} \left( {{\rm {\bf a}}_{.{\kern 1pt}
i}^{ *} } \right), \,\,\,{\mathbb{R}} _{\,{\kern 1pt} ij} =
{\rm{rdet}}_{i} ({\rm {\bf A}}{\rm {\bf A}}^{\ast})_{i.} \left(
{{\rm {\bf a}}_{j.}^{ *} }  \right),\] for $i,j= \overline{ 1,n}$, and ${\rm{ddet}}{ \rm{\bf A}}=\det ({\rm {\bf A}}{\rm {\bf A}}^{\ast})=\det ({\rm {\bf A}}^{\ast}{\rm {\bf A}}).$
\end{theorem}
\subsection{ Determinantal representations of the Moore-Penrose and Drazin inverses over the quaternion skew field}
We shall use the following notations. Let $\alpha : = \left\{
{\alpha _{1} ,\ldots ,\alpha _{k}} \right\} \subseteq {\left\{
{1,\ldots ,m} \right\}}$ and $\beta : = \left\{ {\beta _{1}
,\ldots ,\beta _{k}} \right\} \subseteq {\left\{ {1,\ldots ,n}
\right\}}$ be subsets of the order $1 \le k \le \min {\left\{
{m,n} \right\}}$. By ${\rm {\bf A}}_{\beta} ^{\alpha} $ denote the
submatrix of ${\rm {\bf A}}$ determined by the rows indexed by
$\alpha$ and the columns indexed by $\beta$. Then ${\rm {\bf
A}}{\kern 1pt}_{\alpha} ^{\alpha}$ denotes the principal submatrix
determined by the rows and columns indexed by $\alpha$.
 If ${\rm {\bf A}} \in {\rm
M}\left( {n,{\rm {\mathbb{H}}}} \right)$ is Hermitian, then by
${\left| {{\rm {\bf A}}_{\alpha} ^{\alpha} } \right|}$ denote the
corresponding principal minor of $\det {\rm {\bf A}}$.
 For $1 \leq k\leq n$, the collection of strictly
increasing sequences of $k$ integers chosen from $\left\{
{1,\ldots ,n} \right\}$ is denoted by $\textsl{L}_{ k,
n}: = {\left\{ {\,\alpha :\alpha = \left( {\alpha _{1} ,\ldots
,\alpha _{k}} \right),\,{\kern 1pt} 1 \le \alpha _{1} \le \ldots
\le \alpha _{k} \le n} \right\}}$.  For fixed $i \in \alpha $ and $j \in
\beta $, let $I_{r,\,m} {\left\{ {i} \right\}}: = {\left\{
{\,\alpha :\alpha \in L_{r,m} ,i \in \alpha}  \right\}}{\rm ,}
\quad J_{r,\,n} {\left\{ {j} \right\}}: = {\left\{ {\,\beta :\beta
\in L_{r,n} ,j \in \beta}  \right\}}$.

 Denote by ${\rm {\bf a}}_{.j}^{*} $ and ${\rm {\bf
a}}_{i.}^{*} $ the $j$-th column  and the $i$-th row of  ${\rm
{\bf A}}^{*} $ and by ${\rm {\bf a}}_{.j}^{(m)} $ and ${\rm {\bf
a}}_{i.}^{(m)} $ the $j$-th column  and the $i$-th row of  ${\rm
{\bf A}}^{m} $, respectively.

The following theorem give determinantal representations of the Moore-Penrose inverse over the quaternion skew field $\mathbb{H}$.
\begin{theorem} \cite{ky3}\label{theor:det_repr_MP}
If ${\rm {\bf A}} \in {\rm {\mathbb{H}}}_{r}^{m\times n} $, then
the Moore-Penrose inverse  ${\rm {\bf A}}^{ +} = \left( {a_{ij}^{
+} } \right) \in {\rm {\mathbb{H}}}_{}^{n\times m} $ possess the
following determinantal representations:
\begin{equation}
\label{eq:det_repr_A*A}
 a_{ij}^{ +}  = {\frac{{{\sum\limits_{\beta
\in J_{r,\,n} {\left\{ {i} \right\}}} {{\rm{cdet}} _{i} \left(
{\left( {{\rm {\bf A}}^{ *} {\rm {\bf A}}} \right)_{\,. \,i}
\left( {{\rm {\bf a}}_{.j}^{ *} }  \right)} \right){\kern 1pt}
{\kern 1pt} _{\beta} ^{\beta} } } }}{{{\sum\limits_{\beta \in
J_{r,\,\,n}} {{\left| {\left( {{\rm {\bf A}}^{ *} {\rm {\bf A}}}
\right){\kern 1pt} _{\beta} ^{\beta} }  \right|}}} }}},
\end{equation}
or
\begin{equation}
\label{eq:det_repr_AA*} a_{ij}^{ +}  =
{\frac{{{\sum\limits_{\alpha \in I_{r,m} {\left\{ {j} \right\}}}
{{\rm{rdet}} _{j} \left( {({\rm {\bf A}}{\rm {\bf A}}^{ *}
)_{j\,.\,} ({\rm {\bf a}}_{i.\,}^{ *} )} \right)\,_{\alpha}
^{\alpha} } }}}{{{\sum\limits_{\alpha \in I_{r,\,m}}  {{\left|
{\left( {{\rm {\bf A}}{\rm {\bf A}}^{ *} } \right){\kern 1pt}
_{\alpha} ^{\alpha} } \right|}}} }}}.
\end{equation}
for all $i = \overline{1, n} $, $j =\overline{1, m} $.
\end{theorem}
\begin{proposition}\label{theor:repr_Dr} \cite{ca1} If $Ind({\bf A}) = k$, then
${\bf A}^{D}={\bf A}^{k}({\bf A}^{2k+1})^{+}{\bf A}^{k}
$.
\end{proposition}
Denote  by $\hat{{\rm {\bf a}}}_{.s}$ and $\check{{{\rm {\bf a}}}}_{t.}$ the $s$-th column of $
({\bf A}^{ 2k+1})^{*} {\bf A}^{k}=:\hat{{\rm {\bf A}}}=
(\hat{a}_{ij})\in {\mathbb{H}}^{n\times n}$
and  the $t$-th row of $
{\bf A}^{k}({\bf A}^{ 2k+1})^{*} =:\check{{\rm {\bf A}}}=
(\check{a}_{ij})\in {\mathbb{H}}^{n\times n}$, respectively  for all $s,t=\overline
{1,n}$.
Using  the determinantal representations of the Moore-Penrose inverse (\ref{eq:det_repr_A*A}) and (\ref{eq:det_repr_AA*}), and  Proposition \ref{theor:repr_Dr}  the following determinantal representations of the Drazin inverse for an arbitrary square matrix over $\mathbb{H}$ have been obtained in \cite{ky4}.
\begin{theorem} \cite{ky4}\label{theor:det_rep_draz} If ${\rm {\bf A}} \in {\rm M}\left( {n, {\mathbb{H}}}\right)$  with
$ Ind{\kern 1pt} {\rm {\bf A}}=k$ and $\rank{\rm {\bf A}}^{k+1} =
\rank{\rm {\bf A}}^{k} = r$, then the Drazin inverse  ${\rm {\bf A}}^{ D} $ possess the
determinantal representations,
\begin{equation}
\label{eq:cdet_draz} a_{ij} ^{D}=
 {\frac{{ \sum\limits_{t = 1}^{n} {a}_{it}^{(k)}   {\sum\limits_{\beta \in J_{r,\,n} {\left\{ {t}
\right\}}} {{\rm{cdet}} _{t} \left( {\left({\bf A}^{ 2k+1} \right)^{*}\left({\bf A}^{ 2k+1} \right)_{. \,t} \left( \hat{{\rm {\bf a}}}_{.\,j}
\right)} \right){\kern 1pt}  _{\beta} ^{\beta} } }
}}{{{\sum\limits_{\beta \in J_{r,\,n}} {{\left| {\left({\bf A}^{ 2k+1} \right)^{*}\left({\bf A}^{ 2k+1} \right){\kern 1pt} _{\beta} ^{\beta}
}  \right|}}} }}}
\end{equation}
and
\begin{equation}
\label{eq:rdet_draz} a_{ij} ^{D}=
{\frac{\sum\limits_{s = 1}^{n}\left({{\sum\limits_{\alpha \in I_{r,\,n} {\left\{ {s}
\right\}}} {{\rm{rdet}} _{s} \left( {\left( { {\bf A}^{ 2k+1} \left({\bf A}^{ 2k+1} \right)^{*}
} \right)_{\,. s} (\check{{\rm {\bf a}}}_{i\,.})} \right) {\kern 1pt} _{\alpha} ^{\alpha} } }
}\right){a}_{sj}^{(k)}}{{{\sum\limits_{\alpha \in I_{r,\,n}} {{\left| {\left( { {\bf A}^{ 2k+1} \left({\bf A}^{ 2k+1} \right)^{*}
} \right){\kern 1pt} _{\alpha} ^{\alpha}
}  \right|}}} }}}.
\end{equation}

\end{theorem}
In the special case, when ${\rm {\bf A}} \in {\rm M}\left( {n, {\mathbb{H}}}\right)$ is Hermitian, we can obtain  simpler determinantal
representations of the Drazin inverse.
\begin{theorem} \cite{ky4}\label{theor:det_rep_draz_her}
If ${\rm {\bf A}} \in {\rm M}\left( {n, {\mathbb{H}}}\right)$ is Hermitian with
$ Ind{\kern 1pt} {\rm {\bf A}}=k$ and $\rank{\rm {\bf A}}^{k+1} =
\rank{\rm {\bf A}}^{k} = r$, then the Drazin inverse ${\rm {\bf
A}}^{D} = \left( {a_{ij}^{D} } \right) \in {\rm
{\mathbb{H}}}_{}^{n\times n} $ possess the following determinantal
representations,
\begin{equation}
\label{eq:dr_rep_cdet} a_{ij}^{D}  = {\frac{{{\sum\limits_{\beta
\in J_{r,\,n} {\left\{ {i} \right\}}} {{\rm{cdet}} _{i} \left(
{\left( {{\rm {\bf A}}^{k+1}} \right)_{\,. \,i} \left( {{\rm {\bf
a}}_{.j}^{ k} }  \right)} \right){\kern 1pt} {\kern 1pt} _{\beta}
^{\beta} } } }}{{{\sum\limits_{\beta \in J_{r,\,\,n}} {{\left|
{\left( {{\rm {\bf A}}^{k+1}} \right){\kern 1pt} _{\beta} ^{\beta}
}  \right|}}} }}},
\end{equation}
or
\begin{equation}
\label{eq:dr_rep_rdet} a_{ij}^{D}  = {\frac{{{\sum\limits_{\alpha
\in I_{r,n} {\left\{ {j} \right\}}} {{\rm{rdet}} _{j} \left(
{({\rm {\bf A}}^{ k+1} )_{j\,.\,} ({\rm {\bf a}}_{i.\,}^{ (k)} )}
\right)\,_{\alpha} ^{\alpha} } }}}{{{\sum\limits_{\alpha \in
I_{r,\,n}}  {{\left| {\left( {{\rm {\bf A}}^{k+1} } \right){\kern
1pt}  _{\alpha} ^{\alpha} } \right|}}} }}}.
\end{equation}
\end{theorem}

\subsection{Determinantal representations of the W-weighted Drazin inverse}

\begin{definition}
For an arbitrary matrix over the quaternion skew field, ${\bf A}\in  {\mathbb{H}}^{m\times n}$, we denote by

 $\mathcal{R}_{r}({\rm {\bf A}})=\{ {\bf y}\in {\mathbb{H}}^{m} : \,\,{\bf y} = {\bf A}{\bf x},\,\,  {\bf x} \in {\mathbb{H}}^{n}\}$ , the column right space of ${\bf A}$,

 $\mathcal{N}_{r}({\rm {\bf A}})=\{ {\bf y}\in {\mathbb{H}}^{n} : \,\, {\bf A}{\bf x}=0\}$,  the right null space of  ${\bf A}$,

 $\mathcal{R}_{l}({\rm {\bf A}})=\{ {\bf y}\in {\mathbb{H}}^{n} : \,\,{\bf y} = {\bf x}{\bf A},\,\,  {\bf x} \in {\mathbb{H}}^{m}\}$, the column left space of ${\bf A}$,

 $\mathcal{N}_{r}({\rm {\bf A}})=\{ {\bf y}\in {\mathbb{H}}^{m} : \,\, {\bf x}{\bf A}=0\}$,  the left null space of  ${\bf A}$.
\end{definition}

We introduce some mathematical background from the theory of the W-weighted Drazin inverse \cite{wei,ra,wei1} that can be generalized to ${\mathbb{H}}$.
\begin{lemma}\label{lem:prop}
Let ${\bf A}\in {\mathbb{H}}^{m\times n}$ and ${\bf W}\in {\mathbb{H}}^{n\times m}$ with $k= {\rm max}\{Ind({\bf A}{\bf W}), Ind({\bf W}{\bf A})\}$. Than we have:
\[\begin{aligned}
\;\;\;({\rm a})\;\; & {\rm {\bf A}}_{d,{\bf W}}={\bf A}\left(({\bf W}{\bf A})^{D}) \right)^{2}=\left(({\bf A}{\bf W})^{D}) \right)^{2}{\bf A};\\
\;\;\;({\rm b})\;\; & {\rm {\bf A}}_{d,{\bf W}}{\bf W}=({\rm {\bf A}}{\bf W})^{D};\,\,{\bf W}{\rm {\bf A}}_{d,{\bf W}}=({\bf W}{\rm {\bf A}})^{D};\\
\;\;\;({\rm c})\;\; &  {\bf A}_{d,{\bf W}}= \left\{ ({\bf A}{\bf W})^{k}\left[({\bf A}{\bf W})^{2k+1}\right]^{+}({\bf A}{\bf W})^{k}\right\}{\bf W}^{+};\\
\;\;\;& {\bf A}_{d,{\bf W}}= {\bf W}^{+}\left\{ ({\bf W}{\bf A})^{k}\left[({\bf W}{\bf A})^{2k+1}\right]^{+}({\bf W}{\bf A})^{k}\right\};\\
\;\;\;({\rm d})\;\;\; & {\bf W}{\rm {\bf A}}{\bf W}{\rm {\bf A}}_{d,{\bf W}}={\bf P}_{\mathcal{R}_{r}(({\bf W}{\rm {\bf A}})^{k}),\mathcal{N}_{r}(({\bf W}{\rm {\bf A}})^{k})};\,\,{\rm {\bf A}}_{d,{\bf W}}{\bf W}{\rm {\bf A}}{\bf W}={\bf P}_{\mathcal{R}_{l}(({\rm {\bf A}}{\bf W})^{k}),\mathcal{N}_{l}(({\rm {\bf A}}{\bf W})^{k})},\\
\end{aligned}\]
where ${\bf P}_{\mathcal{R}_{r}(({\bf W}{\rm {\bf A}})^{k}),\mathcal{N}_{r}(({\bf W}{\rm {\bf A}})^{k})}$ is the
 projector on $\mathcal{R}_{r}(({\bf W}{\rm {\bf A}})^{k})$ along $\mathcal{N}_{r}(({\bf W}{\rm {\bf A}})^{k})$,
and ${\bf P}_{\mathcal{R}_{l}(({\rm {\bf A}}{\bf W})^{k}),\mathcal{N}_{l}(({\rm {\bf A}}{\bf W})^{k})}$ is the projector
on $\mathcal{R}_{l}(({\rm {\bf A}}{\bf W})^{k})$ along $\mathcal{N}_{l}(({\rm {\bf A}}{\bf W})^{k})$.
\end{lemma}
In particular, the point (a) of Lemma \ref{lem:prop} due to Cline and Greville \cite{cl} is generalized \cite{song_ela} to ${\mathbb{H}}.$ Using this proposition, we have obtained \cite{ky11}  the following determinantal representations W-weighted Drazin inverse.

Denote ${\bf W}{\bf A}=:{\bf U}=
({u}_{ij})\in {\mathbb{H}}^{n\times n}$ and ${\bf A}{\bf W}=:{\bf V}=
({v}_{ij})\in {\mathbb{H}}^{m\times m}$.

Due to Theorem \ref{theor:det_rep_draz}, we denote an entry of the Drazin inverse ${\bf U}^D$ by
\begin{equation}
\label{eq:u_cdet_draz} u_{ij} ^{D,1}=
 {\frac{{ \sum\limits_{t = 1}^{n} {u}_{it}^{(k)}   {\sum\limits_{\beta \in J_{r,\,n} {\left\{ {t}
\right\}}} {{\rm{cdet}} _{t} \left( {\left({\bf U}^{ 2k+1} \right)^{*}\left({\bf U}^{ 2k+1} \right)_{. \,t} \left( \hat{{\rm {\bf u}}}_{.\,j}
\right)} \right){\kern 1pt}  _{\beta} ^{\beta} } }
}}{{{\sum\limits_{\beta \in J_{r,\,n}} {{\left| {\left({\bf U}^{ 2k+1} \right)^{*}\left({\bf U}^{ 2k+1} \right){\kern 1pt} _{\beta} ^{\beta}
}  \right|}}} }}}
\end{equation}
or
\begin{equation}
\label{eq:u_rdet_draz} u_{ij} ^{D,2}=
{\frac{\sum\limits_{s = 1}^{n}\left({{\sum\limits_{\alpha \in I_{r,\,n} {\left\{ {s}
\right\}}} {{\rm{rdet}} _{s} \left( {\left( { {\bf U}^{ 2k+1} \left({\bf U}^{ 2k+1} \right)^{*}
} \right)_{\,. s} (\check{{\rm {\bf u}}}_{i\,.})} \right) {\kern 1pt} _{\alpha} ^{\alpha} } }
}\right){u}_{sj}^{(k)}}{{{\sum\limits_{\alpha \in I_{r,\,n}} {{\left| {\left( { {\bf U}^{ 2k+1} \left({\bf U}^{ 2k+1} \right)^{*}
} \right){\kern 1pt} _{\alpha} ^{\alpha}
}  \right|}}} }}}
\end{equation}
where $\hat{{\rm {\bf u}}}_{.s}$ and $\check{{{\rm {\bf u}}}}_{t.}$ are the $s$-th column of $
({\bf U}^{ 2k+1})^{*} {\bf U}^{k}=:\hat{{\rm {\bf U}}}=
(\hat{u}_{ij})\in {\mathbb{H}}^{n\times n}$
and  the $t$-th row of $
{\bf U}^{k}({\bf U}^{ 2k+1})^{*} =:\check{{\rm {\bf U}}}=
(\check{u}_{ij})\in {\mathbb{H}}^{n\times n}$, respectively  for all $s,t=\overline
{1,n}$, $r=\rank{\rm {\bf U}}^{k+1} =
\rank{\rm {\bf U}}^{k}$.

Then  we have  the following determinantal representations of ${\rm {\bf A}}_{d,{\bf W}}=
({a}_{ij}^{d,{\bf W}})\in {\mathbb{H}}^{m\times n}$,
\begin{equation}
\label{eq:det_rep1_wdraz}
{a}_{ij}^{d,{\bf W}}=\sum\limits_{q = 1}^{n} a_{iq}(u_{qj}^{D})^{(2)}
\end{equation}
where \begin{equation}
\label{eq:det_rep1_wdraz1}(u_{qj}^{D})^{(2)}=\sum\limits_{p = 1}^{n}u_{qp}^{D,l}u_{pj}^{D,f}\end{equation} for all $l,f=\overline
{1,2}$, and $u_{ij} ^{D,1}$ from (\ref{eq:u_cdet_draz}) and $u_{ij} ^{D,2}$ from (\ref{eq:u_rdet_draz}).

Similarly using ${\bf V}=
({v}_{ij})\in {\mathbb{H}}^{m\times m}$,
\begin{equation}
\label{eq:det_rep2_wdraz}
{a}_{ij}^{d,{\bf W}}=\sum\limits_{q = 1}^{m} (v_{iq}^{D})^{(2)}a_{qj}.
\end{equation}
where the first factor is one of the four possible equations
\begin{equation}
\label{eq:det_rep2_wdraz1}(v_{iq}^{D})^{(2)}=\sum\limits_{p = 1}^{m}v_{ip}^{D,l}v_{pq}^{D,f}\end{equation}
for all $l,f=\overline
{1,2}$, and an entry of the Drazin inverse ${\bf V}^D$ is denoting by
\begin{equation}
\label{eq:v_cdet_draz} v_{ij} ^{D,1}=
 {\frac{{ \sum\limits_{t = 1}^{m} {v}_{it}^{(k)}   {\sum\limits_{\beta \in J_{r,\,m} {\left\{ {t}
\right\}}} {{\rm{cdet}} _{t} \left( {\left({\bf V}^{ 2k+1} \right)^{*}\left({\bf V}^{ 2k+1} \right)_{. \,t} \left( \hat{{\rm {\bf v}}}_{.\,j}
\right)} \right){\kern 1pt}  _{\beta} ^{\beta} } }
}}{{{\sum\limits_{\beta \in J_{r,\,m}} {{\left| {\left({\bf V}^{ 2k+1} \right)^{*}\left({\bf V}^{ 2k+1} \right){\kern 1pt} _{\beta} ^{\beta}
}  \right|}}} }}}
\end{equation}
or
\begin{equation}
\label{eq:v_rdet_draz} v_{ij} ^{D,2}=
{\frac{\sum\limits_{s = 1}^{m}\left({{\sum\limits_{\alpha \in I_{r,\,m} {\left\{ {s}
\right\}}} {{\rm{rdet}} _{s} \left( {\left( { {\bf V}^{ 2k+1} \left({\bf V}^{ 2k+1} \right)^{*}
} \right)_{\,. s} (\check{{\rm {\bf v}}}_{i\,.})} \right) {\kern 1pt} _{\alpha} ^{\alpha} } }
}\right){v}_{sj}^{(k)}}{{{\sum\limits_{\alpha \in I_{r,\,m}} {{\left| {\left( { {\bf V}^{ 2k+1} \left({\bf V}^{ 2k+1} \right)^{*}
} \right){\kern 1pt} _{\alpha} ^{\alpha}
}  \right|}}} }}},
\end{equation}
where $\hat{{\rm {\bf v}}}_{.s}$ and $\check{{{\rm {\bf v}}}}_{t.}$ are the $s$-th column of $
({\bf V}^{ 2k+1})^{*} {\bf V}^{k}=:\hat{{\rm {\bf V}}}=
(\hat{v}_{ij})\in {\mathbb{H}}^{m\times m}$
and  the $t$-th row of $
{\bf V}^{k}({\bf V}^{ 2k+1})^{*} =:\check{{\rm {\bf V}}}=
(\check{v}_{ij})\in {\mathbb{H}}^{m\times m}$, respectively  for all $s,t=\overline
{1,m}$, $r=\rank{\rm {\bf V}}^{k+1} =
\rank{\rm {\bf V}}^{k}$.

 The point (c) of Lemma \ref{lem:prop} due to  \cite{ze} has been generalized  to ${\mathbb{H}}$ in \cite{ky4}. Using this proposition, we have obtained the following two determinantal representations of the W-weighted Drazin inverse.
 \begin{theorem}\cite{ky11}\label{theor:det_rep_vwdraz} Let ${\bf A}\in {\mathbb{H}}^{m\times n}$ and ${\bf W}\in {\mathbb{H}}_{r_{1}}^{n\times m}$ with $k=  Ind({\bf A}{\bf W})$ and $r=\rank({\bf A}{\bf W})^{k+1} =
\rank({\bf A}{\bf W})^{k}$.  Then  the W-weighted Drazin inverse of ${\bf A}$ with respect to ${\bf W}$ possesses the determinantal representations
\begin{multline}
\label{eq:det_repr_v wdraz} {a}_{ij}^{d,{\bf W}}=\\
\frac{\sum\limits_{t = 1}^{m}{\sum\limits_{\alpha \in I_{r,\,m} {\left\{ {t}
\right\}}} {{\rm{rdet}} _{t} \left( {\left( { {\bf V}^{ 2k+1} \left({\bf V}^{ 2k+1} \right)^{*}
} \right)_{t. } (\check{{\rm {\bf v}}}_{i\,.})} \right) {\kern 1pt} _{\alpha} ^{\alpha} } }
{\sum\limits_{\alpha \in I_{r_{1},\,n} {\left\{ {j}
\right\}}} {{\rm{rdet}} _{j} \left( {\left(  {\bf W} {\bf W}^{*} \right)_{j\,. } (\check{{\rm {\bf w}}}_{t.})} \right) {\kern 1pt} _{\alpha} ^{\alpha} } }}{{\sum\limits_{\alpha \in I_{r,\,m}} {{\left| {\left( { {\bf V}^{ 2k+1} \left({\bf V}^{ 2k+1} \right)^{*}
} \right){\kern 1pt} _{\alpha} ^{\alpha}
}  \right|}}}{\sum\limits_{\alpha \in I_{r_{1},\,n}} {{\left| {\left( {\bf W} {\bf W}^{*}\right){\kern 1pt} _{\alpha} ^{\alpha}
}  \right|}}}}
\end{multline} and
\begin{multline}
\label{eq:det_repr_u wdraz} {a}_{ij}^{d,{\bf W}}=\\
\frac{\sum\limits_{t = 1}^{n}
{{\sum\limits_{\beta \in J_{r_{1},\,m} {\left\{ {i}
\right\}}} {{\rm{cdet}} _{i} \left( {\left(  {\bf W}^{*} {\bf W} \right)_{. t} (\hat{{\rm {\bf w}}}_{.t})} \right) {\kern 1pt} _{\beta} ^{\beta} } }
}
{{\sum\limits_{\beta \in J_{r,\,n} {\left\{ {t}
\right\}}} {{\rm{cdet}} _{t} \left( {\left( { \left({\bf U}^{ 2k+1} \right)^{*}{\bf U}^{ 2k+1}
} \right)_{.t } (\hat{{\rm {\bf u}}}_{.j})} \right) {\kern 1pt} _{\beta} ^{\beta} } }
}}
{{\sum\limits_{\beta \in J_{r_{1},\,m}} {{\left| {\left( {\bf W}^{*} {\bf W}\right){\kern 1pt} _{\beta} ^{\beta}
}  \right|}}}{\sum\limits_{\beta \in J_{r,\,
n}} {{\left| {\left( { \left({\bf U}^{ 2k+1} \right)^{*}{\bf U}^{ 2k+1}
} \right){\kern 1pt} _{\beta} ^{\beta}
}  \right|}}}}
\end{multline}
 where
$\check{{\rm {\bf V}}}={\bf V}^{k}({\bf V}^{ 2k+1})^{*}$, $
\check{{\rm {\bf W}}}={\bf V}^{k}{\bf W}^{*}$, and
$\hat{{\rm {\bf U}}}=({\bf U}^{ 2k+1})^{*}{\bf U}^{k}$, $
\hat{{\rm {\bf W}}}={\bf W}^{*}{\bf U}^{k}$.
\end{theorem}
In the special cases,  when
 ${\bf A}{\bf W}\in {\mathbb{H}}^{m\times m}$ and ${\bf W}{\bf A}\in {\mathbb{H}}^{n\times n}$ are Hermitian, we can obtain simpler determinantal representations of the W-weighted Drazin inverse.
\begin{theorem}\cite{ky11}\label{theor:det_rep_wdraz1}
If ${\rm {\bf A}} \in  {\mathbb{H}}^{m\times n}$, ${\bf W}\in {\mathbb{H}}^{n\times m}$, and ${\bf A}{\bf W}\in {\mathbb{H}}^{m\times m}$ is Hermitian with $k= {\rm max}\{Ind({\bf A}{\bf W}), Ind({\bf W}{\bf A})\}$ and $\rank({\bf A}{\bf W})^{k+1} =
\rank({\bf A}{\bf W})^{k} = r$, then the W-weighted Drazin inverse ${\rm {\bf
A}}_{d,W} = \left( {a_{ij}^{d,W} } \right) \in {\rm
{\mathbb{H}}}^{m\times n} $ with respect to ${\bf W}$ possess the following determinantal
representations:
\begin{equation}
\label{eq:dr_rep_wcdet} a_{ij}^{d,W}  = {\frac{{{\sum\limits_{\beta
\in J_{r,\,m} {\left\{ {i} \right\}}} {{\rm{cdet}} _{i} \left(
{\left( {{\rm {\bf A}}{\bf W}} \right)^{k+2}_{\,. \,i} \left( {{\rm {\bf
\bar{v}}}_{.j} }  \right)} \right){\kern 1pt} {\kern 1pt} _{\beta}
^{\beta} } } }}{{{\sum\limits_{\beta \in J_{r,\,\,m}} {{\left|
{\left( { {\bf A}{\bf W}} \right)^{k+2}{\kern 1pt} _{\beta} ^{\beta}
}  \right|}}} }}},
\end{equation}
where ${\rm {\bf \bar{v}}}_{.j} $ is the $j$th column of  ${\rm
{\bf \bar{V}}}=({\bf A}{\bf W})^{k}{\bf A} $ for all $j=\overline{1,m}$.
\end{theorem}
\begin{theorem}\cite{ky11}\label{theor:det_rep_wdraz2}
If ${\rm {\bf A}} \in  {\mathbb{H}}^{m\times n}$, ${\bf W}\in {\mathbb{H}}^{n\times m}$, and ${\bf W}{\bf A}\in {\mathbb{H}}^{n\times n}$ is Hermitian with $k= {\rm max}\{Ind({\bf A}{\bf W}), Ind({\bf W}{\bf A})\}$ and $\rank({\bf W}{\bf A})^{k+1} =
\rank({\bf W}{\bf A})^{k} = r$, then the W-weighted Drazin inverse ${\rm {\bf
A}}_{d,W} = \left( {a_{ij}^{d,W} } \right) \in {\rm
{\mathbb{H}}}^{m\times n} $ with respect to ${\bf W}$ possess the following determinantal
representations:
\begin{equation}
\label{eq:dr_rep_wrdet} a_{ij}^{d,W}  = {\frac{{{\sum\limits_{\alpha
\in I_{r,n} {\left\{ {j} \right\}}} {{\rm{rdet}} _{j} \left(
{({\bf W}{\rm {\bf A}} )^{ k+2}_{j\,.\,} ({\rm {\bf \bar{u}}}_{i.\,}^{ (k)} )}
\right)\,_{\alpha} ^{\alpha} } }}}{{{\sum\limits_{\alpha \in
I_{r,\,n}}  {{\left| {\left({\bf W} {{\rm {\bf A}} } \right)^{k+2}{\kern
1pt}  _{\alpha} ^{\alpha} } \right|}}} }}}.
\end{equation}
where  ${\rm {\bf
\bar{u}}}_{i.} $ is  the $i$th row of  ${\rm
{\bf \bar{U}}}={\bf A}({\bf W}{\bf A})^{k} $ for all $i=\overline{1,
n}$.
\end{theorem}

\section{Cramer's rule for the W-weighted Drazin inverse solution}
\subsection{Background of the problem}
In \cite{wei} Wei has established a Cramer's rule for solving of a general restricted  equation
 \begin{equation}\label{eq:rest_eq}
 {\bf W}{\bf A}{\bf W}{\bf x}={\bf b},\,\,\,{\bf x}\in \mathcal{R}\left[({\bf A}{\bf W})^{k_{1}}\right],
 \end{equation}
  where ${\bf A}\in  {\mathbb{C}}^{m\times n}$, ${\bf W}\in  {\mathbb{C}}^{n\times m}$ with $ Ind\,( {\rm {\bf A}}{\bf W})=k_{1}$,  $ Ind\,( {\bf W}{\rm {\bf A}})=k_{2}$ and $\rank\,( {\rm {\bf A}}{\bf W})^{k_{1}} =r_{1}$, $\rank\,({\bf W} {\rm {\bf A}})^{k_{2}} =r_{2}$.
He proofed that the restricted
matrix equation (\ref{eq:rest_eq}) has a unique solution,
$
 {\bf x}={\bf A}_{d,W}{\bf b}$, and presented its Cramer's rule as follows,
\begin{equation}\label{eq:rest_eq_cr_rul}
 x_j=\det\begin{pmatrix} {\bf W}{\bf A}{\bf W}(j\longrightarrow {\bf b}) & {\bf U}_{1} \\ {\bf V}_{1}(j\longrightarrow 0) & 0 \end{pmatrix}\bigg/ \det\begin{pmatrix} {\bf W}{\bf A}{\bf W} & {\bf U}_{1} \\ {\bf V}_{1} & 0 \end{pmatrix},
  \end{equation}
where ${\bf U}_{1}\in  {\mathbb{C}}^{n\times n-r_{2}}_{n-r_{2}}$,  ${\bf V}^{*}_{1}\in  {\mathbb{C}}^{m\times m-r_{1}}_{m-r_{1}}$ are  matrices whose columns form bases for $\mathcal{N}(({\bf W} {\rm {\bf A}})^{k_{2}})$ and $\mathcal{N}(({\rm {\bf A}}{\bf W} )^{k_{1}^{*}})$, respectively.

Recently, within the framework of a theory
of the column and row determinants Song \cite{song_ela} has considered the characterization of the W-weighted Drazin
inverse over the quaternion skew and presented a Cramer's rule of the restricted matrix equation,
\begin{gather}\label{eq:axb}
{\bf W}_{1}{\bf A}{\bf W}_{1}{\bf X}{\bf W}_{2}{\bf B}{\bf W}_{2}={\bf D},\\
\mathcal{R}_{r}({\rm {\bf X}})\subset\mathcal{R}_{r}\left(({\bf A}{\bf W}_{1})^{k_{1}}\right), \,\,
\mathcal{N}_{r}({\rm {\bf X}})\supset\mathcal{N}_{r}\left(({\bf W}_{2}{\bf B})^{k_{2}}\right),\\ \notag
\label{eq:axbr}
\mathcal{R}_{l}({\rm {\bf X}})\subset\mathcal{R}_{l}\left(({\bf B}{\bf W}_{2})^{k_{2}}\right), \,\,
\mathcal{N}_{l}({\rm {\bf X}})\supset\mathcal{N}_{l}\left(({\bf W}_{1}{\bf A})^{k_{1}}\right),
\end{gather}
where ${\bf A}\in  {\mathbb{H}}^{m\times n}$, ${\bf W}_{1}\in  {\mathbb{H}}^{n\times m}$, ${\bf B}\in  {\mathbb{H}}^{p\times q}$,  ${\bf W}_{2}\in  {\mathbb{H}}^{q\times p}$, and ${\bf D}\in  {\mathbb{H}}^{n\times p}$ with $k_{1}=\max\,\{Ind({\bf A}{\bf W}_{1}),Ind\,({\bf W}_{1}{\bf A})\}$,
$k_{2}=\max\,\{Ind({\bf B}{\bf W}_{2}),Ind\,({\bf W}_{2}{\bf B})\}$, and $\rank\,( {\rm {\bf A}}{\bf W}_{1})^{k_{1}} =s_{1}$, $\rank\,( {\rm {\bf B}}{\bf W}_{2})^{k_{2}} =s_{2}$.

He proofed that if
 \[{\bf D}\in \mathcal{R}_{r}\left(({\bf W}_{1}{\bf A})^{k_{1}}, ({\bf W}_{2}{\bf B})^{k_{2}}\right),\,\, {\bf D}\in \mathcal{R}_{l}\left(({\bf A}{\bf W}_{1})^{k_{1}}, ({\bf B}{\bf W}_{2})^{k_{2}}\right)\]
  and there exist auxiliary matrices of full column rank, ${\bf L}_{1}\in  {\mathbb{H}}^{n\times n-s_{1}}_{n-s_{1}}$, ${\bf M}_{1}^{*}\in  {\mathbb{H}}^{m\times m-s_{1}}_{m-s_{1}}$, ${\bf L}_{2}\in  {\mathbb{H}}^{q\times q-s_{2}}_{q-s_{2}}$, ${\bf M}_{2}^{*}\in  {\mathbb{H}}^{p\times p-s_{2}}_{p-s_{2}}$ with additional terms of their ranges and null spaces, then the restricted
matrix equation (\ref{eq:axb}) has a unique solution, \[ {\bf X}={\bf A}_{d,{\bf W}_{1}}{\bf D}{\bf B}_{d,{\bf W}_{2}}.\]
 Using auxiliary matrices, ${\bf L}_{1}$, ${\bf M}_{1}$, ${\bf L}_{2}$, ${\bf M}_{2}$, Song presented its Cramer's rule  by analogy to (\ref{eq:rest_eq_cr_rul}).

 In this paper we have avoided such approach and have obtained explicit formulas for determinantal representations of the W-weighted Drazin inverse solutions of matrix equations by using only given matrices.

\subsection{A Cramer's rule for the W-weighted Drazin inverse solutions of some matrix equations}
Consider  the following restricted matrix equation,
\begin{gather}\label{eq:ax}
{\bf W}{\bf A}{\bf W}{\bf X}={\bf D},\\
\label{eq:ax_restr}
\mathcal{R}_{r}({\rm {\bf X}})\subset\mathcal{R}_{r}\left(({\bf A}{\bf W})^{k}\right), \,\,
\mathcal{N}_{l}({\rm {\bf X}})\supset\mathcal{N}_{l}\left(({\bf W}{\bf A})^{k}\right),
\end{gather}
where ${\bf A}\in  {\mathbb{H}}^{m\times n}$, ${\bf W}\in  {\mathbb{H}}^{n\times m}_{r_{1}}$ with $k=\max\,\{Ind({\bf A}{\bf W}),Ind\,({\bf W}{\bf A})\}$, and ${\bf D}\in  {\mathbb{H}}^{n\times p}$ .
\begin{theorem}\label{theor:cr_ax}
If ${\rm {\bf D}}\subset\mathcal{R}_{r}\left(({\bf A}{\bf W})^{k}\right)$ and ${\bf D}\supset\mathcal{N}_{l}\left(({\bf W}{\bf A})^{k}\right)$, then the restricted
matrix equation (\ref{eq:ax}) has a unique solution,
\begin{equation}\label{eq:cr_ax}
{\bf X}={\bf A}_{d,W}{\bf D},
\end{equation}
 which possess the following determinantal representations for all $i=\overline{1, m}$,  $j=\overline{1, p}$,

i) \begin{multline}
\label{eq:cr_rul_ax1} {x}_{ij}=\\
\frac{\sum\limits_{t = 1}^{n}
{{\sum\limits_{\beta \in J_{r_{1},\,m} {\left\{ {i}
\right\}}} {{\rm{cdet}} _{i} \left( {\left(  {\bf W}^{*} {\bf W} \right)_{. t} (\hat{{\rm {\bf w}}}_{.t})} \right) {\kern 1pt} _{\beta} ^{\beta} } }
}
{{\sum\limits_{\beta \in J_{r,\,n} {\left\{ {t}
\right\}}} {{\rm{cdet}} _{t} \left( {\left( { \left({\bf U}^{ 2k+1} \right)^{*}{\bf U}^{ 2k+1}
} \right)_{.t } (\hat{{\rm {\bf d}}}_{.j})} \right) {\kern 1pt} _{\beta} ^{\beta} } }
}}
{{\sum\limits_{\beta \in J_{r_{1},\,m}} {{\left| {\left( {\bf W}^{*} {\bf W}\right){\kern 1pt} _{\beta} ^{\beta}
}  \right|}}}{\sum\limits_{\beta \in J_{r,\,
n}} {{\left| {\left( { \left({\bf U}^{ 2k+1} \right)^{*}{\bf U}^{ 2k+1}
} \right){\kern 1pt} _{\beta} ^{\beta}
}  \right|}}}}
\end{multline}
where $\hat{{\rm {\bf d}}}_{.j}$ is the $j$-th column of $\hat{{\rm {\bf D}}}= \hat{{\rm {\bf U}}}{\bf D}=({\bf U}^{2k+1})^{*}{\bf U}^{k}{\bf D}$, ${\bf U}={\bf W}{\bf A}$, $
\hat{{\rm {\bf W}}}={\bf W}^{*}{\bf U}^{k}$, and $r=\rank ({\bf W}{\bf A})^{k+1}=\rank ({\bf W}{\bf A})^{k}$.

   ii)\begin{equation}
\label{eq:cr_rul_ax2}{x}_{ij}=\sum\limits_{q = 1}^{m} (v_{iq}^{D})^{(2)}{r}_{qj},\end{equation} where $(v_{iq}^{D})^{(2)}$ can be obtained by (\ref{eq:det_rep2_wdraz1}) and ${\bf A}{\bf D}={\bf R}=({r}_{qj})\in  {\mathbb{H}}^{m\times p}$.

iii) If ${\bf A}{\bf W}\in {\mathbb{H}}^{m\times m}$ is Hermitian, then

\begin{equation}
\label{eq:cr_rul_ax3}{x}_{ij}  = {\frac{{{\sum\limits_{\beta
\in J_{r,\,m} {\left\{ {i} \right\}}} {{\rm{cdet}} _{i} \left(
{\left( {{\rm {\bf A}}{\bf W}} \right)^{k+2}_{\,. \,i} \left( {{\rm {\bf
{f}}}_{.j} }  \right)} \right){\kern 1pt} {\kern 1pt} _{\beta}
^{\beta} } } }}{{{\sum\limits_{\beta \in J_{r,\,\,m}} {{\left|
{\left( { {\bf A}{\bf W}} \right)^{k+2}{\kern 1pt} _{\beta} ^{\beta}
}  \right|}}} }}},
\end{equation}
where ${\rm {\bf {f}}}_{.j} $ is the $j$-th column of  ${{\bf F}}={\rm
{\bf \bar{V}}}{\bf D}=({\bf A}{\bf W})^{k}{\bf A}{\bf D} $.

\end{theorem}
\emph{Proof.} The proof contains two parts. We first shall establish that the unique solution
of (\ref{eq:ax}) can be represented as (\ref{eq:cr_ax}). By the definition of the right range, we have
\[{\bf D}=({\bf W}{\bf A})^{k}{\bf Y}\]
for some matrix ${\bf Y}\in  {\mathbb{H}}^{n\times p}$. It follows that
\[ \mathcal{R}_{r}({\rm {\bf D}})\subset\mathcal{R}_{r}\left(({\bf W}{\bf A})^{k}\right).\]
Then by Lemma \ref{lem:prop} (d),
\[{\bf W}{\bf A}{\bf W}{\bf A}_{d,W}{\bf D}={\bf P}_{\mathcal{R}_{r}(({\bf W}{\rm {\bf A}})^{k}),\mathcal{N}_{r}(({\bf W}{\rm {\bf A}})^{k})}{\bf D}={\bf D}.\]
It means that (\ref{eq:cr_ax}) is a solution of (\ref{eq:ax}) and satisfies the restricted conditions (\ref{eq:ax_restr}).

Now we prove the uniqueness of (\ref{eq:cr_ax}). Let ${\bf X}_{0}$ is a solution of (\ref{eq:ax}). Then it satisfies the restricted conditions (\ref{eq:ax_restr}), and
\[{\bf A}_{d,W}{\bf D}= {\bf A}_{d,W}{\bf W}{\bf A}{\bf W}{\bf X}_{0}={\bf P}_{\mathcal{R}_{r}(({\bf W}{\rm {\bf A}})^{k}),\mathcal{N}_{r}(({\bf W}{\rm {\bf A}})^{k})}{\bf X}_{0}={\bf X}_{0}.\]

To derive a Cramer's rule (\ref{eq:cr_rul_ax1}), we use  the  determinantal representation (\ref{eq:det_repr_u wdraz}) for ${\bf A}_{d,W}$. Then
\begin{multline}\label{eq:cr_rul_ax11}
x_{ij}=\sum\limits_{s = 1}^{p}{a}_{is}^{d,{\bf W}}d_{sj}=\\\sum\limits_{s = 1}^{p} \left[\frac{\sum\limits_{t = 1}^{n}
{{\sum\limits_{\beta \in J_{r_{1},\,m} {\left\{ {i}
\right\}}} {{\rm{cdet}} _{i} {\left(  {\bf W}^{*} {\bf W} \right)_{. t} (\hat{{\rm {\bf w}}}_{.t})}  _{\beta} ^{\beta} } }
}
{{\sum\limits_{\beta \in J_{r,\,n} {\left\{ {t}
\right\}}} {{\rm{cdet}} _{t}  {\left( { \left({\bf U}^{ 2k+1} \right)^{*}{\bf U}^{ 2k+1}
} \right)_{.t } (\hat{{\rm {\bf u}}}_{.s})} _{\beta} ^{\beta} } }
}}
{{\sum\limits_{\beta \in J_{r_{1},\,m}} {{\left| {\left( {\bf W}^{*} {\bf W}\right){\kern 1pt} _{\beta} ^{\beta}
}  \right|}}}{\sum\limits_{\beta \in J_{r,\,
n}} {{\left| {\left( { \left({\bf U}^{ 2k+1} \right)^{*}{\bf U}^{ 2k+1}
} \right){\kern 1pt} _{\beta} ^{\beta}
}  \right|}}}}\right] d_{sj}
\end{multline}
Denote $\hat{{\rm {\bf D}}}= \hat{{\rm {\bf U}}}{\bf D}=({\bf U}^{2k+1})^{*}{\bf U}^{k}{\bf D}$, where $\hat{{\rm {\bf D}}}=\left(\hat{d}_{sj}\right)\in {\mathbb{H}}^{n\times p}$. Since \[\sum\limits_{s = 1}^{p}\hat{{\rm {\bf u}}}_{.s}d_{sj}=\hat{{\bf d}}_{.j},\]

where $\hat{{\bf d}}_{.j}$ is the $j$-th column of $\hat{{\rm {\bf D}}}$, then (\ref{eq:cr_rul_ax1}) follows from (\ref{eq:cr_rul_ax11}).

Similarly, we derive the analogs of Cramer's rule (\ref{eq:cr_rul_ax2}) and (\ref{eq:cr_rul_ax3}) by using the determinantal representations for the W-weighted Drazin inverse (\ref{eq:det_rep2_wdraz}), and (\ref{eq:dr_rep_wcdet}), respectively.
$\Box$
\begin{remark}  In the complex case, i.e. ${\bf A}\in  {\mathbb{C}}^{m\times n}$, ${\bf W}\in  {\mathbb{C}}^{n\times m}_{r_{1}}$, and ${\bf D}\in  {\mathbb{C}}^{n\times p}$, we substitute usual determinants for all corresponding row and column determinants in (\ref{eq:cr_rul_ax1}), (\ref{eq:cr_rul_ax2}), and (\ref{eq:cr_rul_ax3}).

Note that in the case iii), the condition ${\bf A}{\bf W}\in {\mathbb{C}}^{m\times m}$ be Hermitian is not necessary, then in the complex case (\ref{eq:cr_rul_ax3})  will have the form
 \begin{equation*}
{x}_{ij}  = {\frac{{{\sum\limits_{\beta
\in J_{r,\,m} {\left\{ {i} \right\}}} \left|{ \left(
{\left( {{\rm {\bf A}}{\bf W}} \right)^{k+2}_{\,. \,i} \left( {{\rm {\bf
{f}}}_{.j} }  \right)} \right){\kern 1pt} {\kern 1pt} _{\beta}
^{\beta} } \right| }}}{{{\sum\limits_{\beta \in J_{r,\,\,m}} {{\left|
{\left( { {\bf A}{\bf W}} \right)^{k+2}{\kern 1pt} _{\beta} ^{\beta}
}  \right|}}} }}},
\end{equation*}
where ${\rm {\bf {f}}}_{.j} $ is the $j$-th column of  ${{\bf F}}={\rm
{\bf \bar{V}}}{\bf D}=({\bf A}{\bf W})^{k}{\bf A}{\bf D} $.
\end{remark}

Now, consider  the following restricted matrix equation,
\begin{equation}\label{eq:xa}
{\bf X}{\bf W}{\bf A}{\bf W}={\bf D},\\
\end{equation}
\begin{equation}\label{eq:xa_restr}
\mathcal{R}_{l}({\rm {\bf X}})\subset\mathcal{R}_{l}\left(({\bf A}{\bf W})^{k}\right), \,\,
\mathcal{N}_{r}({\rm {\bf X}})\supset\mathcal{N}_{r}\left(({\bf W}{\bf A})^{k}\right),
\end{equation}
where ${\bf A}\in  {\mathbb{H}}^{m\times n}$, ${\bf W}\in  {\mathbb{H}}^{n\times m}_{r_{1}}$ with $k=\max\,\{Ind({\bf A}{\bf W}),Ind\,({\bf W}{\bf A})\}$, and ${\bf D}\in  {\mathbb{H}}^{q\times m}$.
\begin{theorem}\label{theor:cr_xa}
If ${\rm {\bf D}}\subset\mathcal{R}_{l}\left(({\bf A}{\bf W})^{k}\right)$ and ${\bf D}\supset\mathcal{N}_{r}\left(({\bf W}{\bf A})^{k}\right)$, then the restricted
matrix equation (\ref{eq:xa}) has a unique solution,
\begin{equation}\label{eq:cr_xa}
{\bf X}={\bf D}{\bf A}_{d,W},
\end{equation}
 which possess the following determinantal representations for  $i=\overline{1, q}$,  $j=\overline{1, n}$,

i)\begin{multline}
\label{eq:cr_rul_xa1} {x}_{ij}=\\
\frac{\sum\limits_{l = 1}^{m}{\sum\limits_{\alpha \in I_{r,\,m} {\left\{ {t}
\right\}}} {{\rm{rdet}} _{l} \left( {\left( { {\bf V}^{ 2k+1} \left({\bf V}^{ 2k+1} \right)^{*}
} \right)_{l. } (\check{{\rm {\bf d}}}_{i\,.})} \right) {\kern 1pt} _{\alpha} ^{\alpha} } }
{\sum\limits_{\alpha \in I_{r_{1},\,n} {\left\{ {j}
\right\}}} {{\rm{rdet}} _{j} \left( {\left(  {\bf W} {\bf W}^{*} \right)_{j\,. } (\check{{\rm {\bf w}}}_{l.})} \right) {\kern 1pt} _{\alpha} ^{\alpha} } }}{{\sum\limits_{\alpha \in I_{r,\,m}} {{\left| {\left( { {\bf V}^{ 2k+1} \left({\bf V}^{ 2k+1} \right)^{*}
} \right){\kern 1pt} _{\alpha} ^{\alpha}
}  \right|}}}{\sum\limits_{\alpha \in I_{r_{1},\,n}} {{\left| {\left( {\bf W} {\bf W}^{*}\right){\kern 1pt} _{\alpha} ^{\alpha}
}  \right|}}}}
\end{multline}
where $\check{{\rm {\bf d}}}_{i.}$ is the $i$-th row of $\check{{\rm {\bf D}}}={\bf D} \check{{\rm {\bf V}}}={\bf D}{\bf V}^{k}({\bf V}^{2k+1})^{*}$, ${\bf V}={\bf A}{\bf W}$ and $r=\rank ({\bf A}{\bf W})^{k+1}=\rank ({\bf A}{\bf W})^{k}$.

   ii)\begin{equation}
\label{eq:cr_rul_xa2}{x}_{ij}=\sum\limits_{t = 1}^{n} l_{it}(u_{tj}^{D})^{(2)},\end{equation} where $(u_{qj}^{D})^{(2)}$ can be obtained by
 (\ref{eq:det_rep1_wdraz1})  and ${\bf D}{\bf A}={\bf L}=({l}_{it})\in  {\mathbb{H}}^{q\times n}$.

iii) If ${\bf A}{\bf W}\in {\mathbb{H}}^{m\times m}$ is Hermitian, then
\begin{equation}
\label{eq:cr_rul_xa3}{x}_{ij}  = {\frac{{{\sum\limits_{\alpha
\in I_{r,n} {\left\{ {j} \right\}}} {{\rm{rdet}} _{j} \left(
{({\bf W}{\rm {\bf A}} )^{ k+2}_{j\,.\,} ( {{\bf{g}}}_{i.\,} )}
\right)\,_{\alpha} ^{\alpha} } }}}{{{\sum\limits_{\alpha \in
I_{r,\,n}}  {{\left| {\left({\bf W} {{\rm {\bf A}} } \right)^{k+2}{\kern
1pt}  _{\alpha} ^{\alpha} } \right|}}} }}}.
\end{equation}
where  ${{\bf{g}}}_{i.} $ is  the $i$-th row of  ${\rm
{\bf {G}}}={\bf D}{\bf A}({\bf W}{\bf A})^{k} $ for all $i=\overline{1,
n}$.
\end{theorem}
\emph{Proof.} The proof is similar to the proof of Theorem \ref{theor:cr_ax}.
$\Box$
\begin{remark}  In the complex case, i.e. ${\bf A}\in  {\mathbb{C}}^{m\times n}$, ${\bf W}\in  {\mathbb{C}}^{n\times m}_{r_{1}}$, and ${\bf D}\in  {\mathbb{C}}^{n\times p}$, we substitute usual determinants for all corresponding row and column determinants in (\ref{eq:cr_rul_xa1}), (\ref{eq:cr_rul_xa2}), and (\ref{eq:cr_rul_xa3}). Herein the condition ${\bf W}{\bf A}\in {\mathbb{C}}^{n\times n}$ be Hermitian is not necessary, then in the complex case (\ref{eq:cr_rul_xa3})  can be represented as follows,
\begin{equation*}
{x}_{ij}  = {\frac{{{\sum\limits_{\alpha
\in I_{r,n} {\left\{ {j} \right\}}}\left| { \left(
{({\bf W}{\rm {\bf A}} )^{ k+2}_{j\,.\,} ( {{\bf{g}}}_{i.\,} )}
\right)\,_{\alpha} ^{\alpha} }\right| }}}{{{\sum\limits_{\alpha \in
I_{r,\,n}}  {{\left| {\left({\bf W} {{\rm {\bf A}} } \right)^{k+2}{\kern
1pt}  _{\alpha} ^{\alpha} } \right|}}} }}}.
\end{equation*}
where  ${{\bf{g}}}_{i.} $ is  the $i$-th row of  ${\rm
{\bf {G}}}={\bf D}{\bf A}({\bf W}{\bf A})^{k} $ for all $i=\overline{1,
n}$.

\end{remark}

Now we consider the matrix equation (\ref{eq:axb}) with the constraints (\ref{eq:axbr}). Denote ${\bf A}{\bf D} {\bf B}=: {\tilde{{\bf D}}}=\left({\tilde{d}}_{lf}\right)\in  {\mathbb{H}}^{m\times q}$, and $  {\bar{\bf
V}}{\rm {\bf D}}{\bar{\bf U}}=:{\bar {\bf {D}}}=\left({\bar{d}}_{lf}\right)\in  {\mathbb{H}}^{m\times q}$, where ${\bar{\bf
V}}:=({\bf A}{\bf W}_{1})^{k_{1}}{\bf A}$ , ${\bar{\bf
U}}:={\bf B}({\bf W}_{2}{\bf B})^{k_{2}}$.
\begin{theorem}\label{theor:cr_axb} Suppose ${\bf D}\in  {\mathbb{H}}^{n\times p}$, ${\bf A}\in  {\mathbb{H}}^{m\times n}$, ${\bf W}_{1}\in  {\mathbb{H}}^{n\times m}_{r_{1}}$ with $k_{1}=\max\,\{Ind({\bf A}{\bf W}_{1}),Ind\,({\bf W}_{1}{\bf A})\}$, where  $\rank\,( {\rm {\bf A}}{\bf W}_{1})^{k_{1}} =s_{1}$, and
 ${\bf B}\in  {\mathbb{H}}^{p\times q}$,  ${\bf W}_{2}\in  {\mathbb{H}}^{q\times p}_{r_{2}}$ with
$k_{2}=\max\,\{Ind({\bf B}{\bf W}_{2}),Ind\,({\bf W}_{2}{\bf B})\}$,   $\rank\,( {\rm {\bf B}}{\bf W}_{2})^{k_{2}} =s_{2}$.
If ${\bf D}\in \mathcal{R}_{r}\left(({\bf W}_{1}{\bf A})^{k_{1}}, ({\bf W}_{2}{\bf B})^{k_{2}}\right)$, ${\bf D}\in \mathcal{R}_{l}\left(({\bf A}{\bf W}_{1})^{k_{1}}, ({\bf B}{\bf W}_{2})^{k_{2}}\right)$, then the restricted
matrix equation (\ref{eq:axb}) has a unique solution,
\begin{equation}\label{eq:cr_axb}
{\bf X}={\bf A}_{d,{\bf W}_{1}}{\bf D}{\bf B}_{d,{\bf W}_{2}},
\end{equation}
 which possess the following determinantal representations for all $i=\overline{1, m}$,  $j=\overline{1, q}$.

 i)\begin{equation}\label{eq:cr_rul_axb_2}
 {x}_{ij}=\sum\limits_{l = 1}^{m}\sum\limits_{f = 1}^{q} (v_{il}^{D})^{(2)}{\tilde{d}}_{lf}(u_{fj}^{D})^{(2)},\end{equation}
 where  $(v_{il}^{D})={\bf V}^{D}$ is the Drazin inverse of ${\bf V}={\bf A}{\bf W}_{1}$ and $(v_{il}^{D})^{(2)}$ can be obtained by (\ref{eq:det_rep2_wdraz1}), and
 $(u_{fj}^{D})={\bf U}^{D}$ is the Drazin inverse of ${\bf U}={\bf W}_{2}{\bf B}$ and $(u_{qj}^{D})^{(2)}$ can be obtained by
 (\ref{eq:det_rep1_wdraz1}).

ii) If ${\bf A}{\bf W}_{1}\in {\mathbb{H}}^{m\times m}$ and ${\bf W}_{2}{\bf B}\in {\mathbb{H}}^{q\times q}$ are Hermitian, then
\begin{equation}\label{eq:d^B}
x_{ij} = {\frac{{\sum\limits_{\beta \in
J_{s_{1},\,m} {\left\{ {i} \right\}}} {{\rm{cdet}} _{i} \left(
{\left( {{\rm {\bf A}}{\bf W}_{1}} \right)^{k_{1}+2}_{\,. \,i} \left( {{{\rm {\bf d}}}_{.\,j}^{{\rm {\bf B}}}}
\right)} \right){\kern 1pt} _{\beta} ^{\beta} } }}{{{{\sum\limits_{\beta \in J_{s_{1},\,m}} {{\left|
{\left( { {\bf A}{\bf W}_{1}} \right)^{k_{1}+2}{\kern 1pt} _{\beta} ^{\beta}
}  \right|}}} }}{{\sum\limits_{\alpha \in
I_{s_{2},\,q}}  {{\left| {\left({\bf W}_{2} {{\rm {\bf B}} } \right)^{k_{2}+2}{\kern
1pt}  _{\alpha} ^{\alpha} } \right|}}} }}},
\end{equation}
or
\begin{equation}\label{eq:d^A}
 x_{ij}={\frac{\sum\limits_{\alpha \in I_{s_{2},q} {\left\{ {j}
\right\}}} {{\rm{rdet}} _{j} \left( {({\bf W}_{2}{\rm {\bf B}} )^{ k_{2}+2}_{j\,.\,} ({\rm {\bf
d}}^{{\rm {\bf A}}}_{i\,.} )} \right)\,_{\alpha} ^{\alpha}
} }{{{{\sum\limits_{\beta \in J_{s_{1},\,m}} {{\left|
{\left( { {\bf A}{\bf W}_{1}} \right)^{k_{1}+2}{\kern 1pt} _{\beta} ^{\beta}
}  \right|}}} }}{{\sum\limits_{\alpha \in
I_{s_{2},\,q}}  {{\left| {\left({\bf W}_{2} {{\rm {\bf B}} } \right)^{k_{2}+2}{\kern
1pt}  _{\alpha} ^{\alpha} } \right|}}} }}},
\end{equation}
where
\begin{equation} \label{eq:def_d^B_m}
   {{{\rm {\bf d}}}_{.\,j}^{{\rm {\bf B}}}}=\left(
{\sum\limits_{\alpha \in I_{s_{2},q} {\left\{ {j}
\right\}}} {{\rm{rdet}} _{j} \left( {({\bf W}_{2}{\rm {\bf B}} )^{ k_{2}+2}_{j\,.\,} ( {\bar {\bf
d}}_{t.} )} \right)_{\alpha} ^{\alpha}
} }
\right)\in{\rm {\mathbb{H}}}^{n \times
1},\,\,\,\,t=\overline{1,n} \end{equation}
\begin{equation} \label{eq:d_A<n}  {{{\rm {\bf d}}}_{i\,.}^{{\rm {\bf A}}}}=\left(
{\sum\limits_{\beta \in
J_{s_{1},\,m} {\left\{ {i} \right\}}} {{\rm{cdet}} _{i} \left(
{\left( {{\rm {\bf A}}{\bf W}_{1}} \right)^{k_{1}+2}_{\,. \,i} \left( {{\bf {\bar d}}}_{.\,l}
\right)} \right){\kern 1pt} _{\beta} ^{\beta} } } \right)\in{\rm {\mathbb{H}}}^{1 \times
q},\,\,\,\,l=\overline{1,q}
\end{equation}
 are the column vector and the row vector, respectively.  ${{\bf {\bar
 d}}_{i.}}$ and
${{\bf {\bar
 d}}_{.j}}$ are the i-th row   and the j-th
column  of ${\bf {\bar{D}}}$ for all $i =\overline{1,n} $, $j =\overline{1,p} $.

\end{theorem}
\emph{Proof.} The existence and  uniqueness of the solution (\ref{eq:cr_axb}) can be proved similar as in (\cite{{song_ela}}, Theorem 5.2).

To establish a Cramer's rule of  (\ref{eq:axb}) we note that we shall not use the determinantal representations (\ref{eq:dr_rep_wcdet}) and (\ref{eq:dr_rep_wcdet}) for  (\ref{eq:cr_axb})   because corresponding determinantal representations of it's solution will be too cumbersome.

To derive a Cramer's rule (\ref{eq:cr_rul_axb_2}) we use the sentence (a) from Lemma \ref{lem:prop}. Then we obtain
\begin{equation}\label{eq:cr_rul_axb_21}
{\bf X}=\left(({\bf A}{\bf W}_{1})^{D}) \right)^{2}{\bf A}{\bf D}
 {\bf B}\left(({\bf W}_{2}{\bf B})^{D}) \right)^{2}.\end{equation}
Denote ${\bf A}{\bf D} {\bf B}=: {\tilde{{\bf D}}}=\left({\tilde{d}}_{lf}\right)\in  {\mathbb{H}}^{m\times q}$,
${\bf V}:={\bf A}{\bf W}_{1}$, and ${\bf U}:={\bf W}_{2}{\bf B}$.
Then the equation  (\ref{eq:cr_rul_axb_21}) will be written component-wise as follows
 \begin{multline*}{x}_{ij}=\sum\limits_{s = 1}^{p}\sum\limits_{t = 1}^{n} (a_{it}^{d,W{1}}){d}_{ts}(b_{sj}^{d,W_{2}})=\\
\sum\limits_{s = 1}^{p}\sum\limits_{t = 1}^{n}\left( \sum\limits_{l = 1}^{m}(v_{il}^{D})^{(2)}a_{lt}\right)d_{ts} \left(\sum\limits_{f = 1}^{q} b_{sf}(u_{fj}^{D})^{(2)}\right)
\end{multline*}
 By changing the order of summation, from here  it follows (\ref{eq:cr_rul_axb_2}).

ii) If ${\rm {\bf A}} \in {\rm {\mathbb{H}}}_{r_{1}}^{m\times n}
$, ${\rm {\bf B}} \in {\rm {\mathbb{H}}}_{r_{2}}^{p\times q} $ and
${\bf A}{\bf W}_{1}\in {\mathbb{H}}^{m\times m}$ and ${\bf W}_{2}{\bf B}\in {\mathbb{H}}^{q\times q}$ are Hermitian, then by Theorems
\ref{theor:det_rep_wdraz1} and \ref{theor:det_rep_wdraz2} the W-weighted Drazin inverses ${\rm {\bf
A}}_{d,W_{1}} = \left( {a_{ij}^{d,W_{1}} } \right) \in {\rm
{\mathbb{H}}}^{m\times n} $ and ${\rm {\bf B}}_{d,W_{2}} = \left(
{b_{ij}^{d,W_{2}} } \right) \in {\rm {\mathbb{H}}}^{q\times p} $ posses
the following determinantal representations respectively,
\begin{equation}
\label{eq:dr_rep_wcdet2} a_{ij}^{d,W_{1}}  = {\frac{{{\sum\limits_{\beta
\in J_{s_{1},\,m} {\left\{ {i} \right\}}} {{\rm{cdet}} _{i} \left(
{\left( {{\rm {\bf A}}{\bf W_{1}}} \right)^{k_{1}+2}_{\,. \,i} \left( {{\rm {\bf
\bar{v}}}_{.j} }  \right)} \right){\kern 1pt} {\kern 1pt} _{\beta}
^{\beta} } } }}{{{\sum\limits_{\beta \in J_{r,\,\,m}} {{\left|
{\left( { {\bf A}{\bf W_{1}}} \right)^{k_{1}+2}{\kern 1pt} _{\beta} ^{\beta}
}  \right|}}} }}},
\end{equation}
where ${\rm {\bf \bar{v}}}_{.j} $ is the $j$-th column of  ${\rm
{\bf \bar{V}}}=({\bf A}{\bf W_{1}})^{k_{1}}{\bf A} $ for all $j=\overline{1,m}$;
\begin{equation}
\label{eq:dr_rep_wrdet2} b_{ij}^{d,W_{2}}  = {\frac{{{\sum\limits_{\alpha
\in I_{s_{2},q} {\left\{ {j} \right\}}} {{\rm{rdet}} _{j} \left(
{({\bf W_{2}}{\rm {\bf B}} )^{ k_{2}+2}_{j\,.\,} ({\rm {\bf \bar{u}}}_{i.\,} )}
\right)\,_{\alpha} ^{\alpha} } }}}{{{\sum\limits_{\alpha \in
I_{s_{2},\,q}}  {{\left| {\left({\bf W_{2}} {{\rm {\bf B}} } \right)^{k_{2}+2}{\kern
1pt}  _{\alpha} ^{\alpha} } \right|}}} }}},
\end{equation}
where  ${\rm {\bf
\bar{u}}}_{i.} $ is  the $i$-th row of  ${\rm
{\bf \bar{U}}}={\bf B}({\bf W_{2}}{\bf B})^{k_{2}} $ for all $i=\overline{1,p
}$.

By component-wise writing (\ref{eq:cr_axb}) we obtain,
  \begin{gather}\label{eq:cr_axb_comp}
  {x}_{ij}=\sum\limits_{s = 1}^{p}\left(\sum\limits_{t = 1}^{n} a_{it}^{d,W{1}}{d}_{ts}\right)\cdot b_{sj}^{d,W_{2}}
   \end{gather}

   Denote  by $\hat{{\rm {\bf d}}_{.s}}$ the $s$-th column of ${
\bar{\bf V}}{\rm {\bf D}}= ({\bf A}{\bf W_{1}})^{k_{1}}{\bf A}{\bf D}=:\hat{{\rm {\bf D}}}=
(\hat{d}_{ij})\in {\mathbb{H}}^{m\times p}$ for all $s=\overline{1,p}$. It follows from ${\sum\limits_{t}  {\bar {\bf v}}_{.\,t}}d_{ts} =\hat{{\rm {\bf d}}_{.\,s}}$ that
\begin{multline}\label{eq:sum_cdet}
\sum\limits_{t = 1}^{n} a_{it}^{d,W{1}}{d}_{ts}=\sum\limits_{t = 1}^{n}
{\frac{{{\sum\limits_{\beta
\in J_{s_{1},\,m} {\left\{ {i} \right\}}} {{\rm{cdet}} _{i} \left(
{\left( {{\rm {\bf A}}{\bf W}_{1}} \right)^{k_{1}+2}_{\,. \,i} \left( {{\rm {\bf
\bar{v}}}_{.t} }  \right)} \right){\kern 1pt} {\kern 1pt} _{\beta}
^{\beta} } } }}{{{\sum\limits_{\beta \in J_{s_{1},\,\,m}} {{\left|
{\left( { {\bf A}{\bf W}_{1}} \right)^{k_{1}+2}{\kern 1pt} _{\beta} ^{\beta}
}  \right|}}} }}}\cdot {d}_{ts} =\\
{\frac{{{\sum\limits_{\beta
\in J_{s_{1},\,m} {\left\{ {i} \right\}}}\sum\limits_{t = 1}^{n} {{\rm{cdet}} _{i} \left(
{\left( {{\rm {\bf A}}{\bf W}_{1}} \right)^{k_{1}+2}_{\,. \,i} \left( {{\rm {\bf
\bar{v}}}_{.t} }  \right)} \right){\kern 1pt}  _{\beta}
^{\beta} } } }\cdot {d}_{ts}   }{{{\sum\limits_{\beta \in J_{s_{1},\,m}} {{\left|
{\left( { {\bf A}{\bf W}_{1}} \right)^{k_{1}+2}{\kern 1pt} _{\beta} ^{\beta}
}  \right|}}} }}} =\\
{\frac{{{\sum\limits_{\beta
\in J_{s_{1},\,m} {\left\{ {i} \right\}}}{{\rm{cdet}} _{i} \left(
{\left( {{\rm {\bf A}}{\bf W}_{1}} \right)^{k_{1}+2}_{\,. \,i} \left( { \bf { \hat d}}_{.\,s}  \right)} \right){\kern 1pt} _{\beta}
^{\beta} } } }  }{{{\sum\limits_{\beta \in J_{s_{1},\,m}} {{\left|
{\left( { {\bf A}{\bf W}_{1}} \right)^{k_{1}+2}{\kern 1pt} _{\beta} ^{\beta}
}  \right|}}} }}}
\end{multline}
   Suppose ${\rm {\bf e}}_{s.}$ and ${\rm {\bf e}}_{.\,s}$ are
respectively the unit row-vector and the unit column-vector whose
components are $0$, except the $s$-th components, which are $1$.
Substituting  (\ref{eq:sum_cdet}) and (\ref{eq:dr_rep_wrdet2}) in
(\ref{eq:cr_axb_comp}), we obtain
\begin{multline*}
x_{ij} =\\
\sum\limits_{s = 1}^{p}
{\frac{{{\sum\limits_{\beta
\in J_{s_{1},\,m} {\left\{ {i} \right\}}}{{\rm{cdet}} _{i} \left(
{\left( {{\rm {\bf A}}{\bf W}_{1}} \right)^{k_{1}+2}_{\,. \,i} \left( { \bf { \hat d}}_{.\,s}  \right)} \right){\kern 1pt} _{\beta}
^{\beta} } } }  }{{{\sum\limits_{\beta \in J_{s_{1},\,m}} {{\left|
{\left( { {\bf A}{\bf W}_{1}} \right)^{k_{1}+2}{\kern 1pt} _{\beta} ^{\beta}
}  \right|}}} }}}
{\frac{{{\sum\limits_{\alpha
\in I_{s_{2},q} {\left\{ {j} \right\}}} {{\rm{rdet}} _{j} \left(
{({\bf W}_{2}{\rm {\bf B}} )^{ k_{2}+2}_{j\,.\,} ({\rm {\bf \bar{u}}}_{s.\,} )}
\right)\,_{\alpha} ^{\alpha} } }}}{{{\sum\limits_{\alpha \in
I_{s_{2},\,q}}  {{\left| {\left({\bf W}_{2} {{\rm {\bf B}} } \right)^{k_{2}+2}{\kern
1pt}  _{\alpha} ^{\alpha} } \right|}}} }}}.
\end{multline*}
Since \begin{equation}\label{eq:prop}{{\bf {\hat
d}}_{.\,s}}=\sum\limits_{t = 1}^{n}{\rm {\bf e}}_{.\,t}\hat{
d_{ts}},\,   {\bf \bar{u}}_{s.\,}=\sum\limits_{l =
1}^{q}{\bar u}_{sl}{\rm {\bf
e}}_{l.},\,\sum\limits_{s=1}^{p}\hat{d_{ts}}{\bar u}_{sl}=\bar {d}_{tl},\end{equation}
then we have
\begin{multline}\label{eq:x_ij}
x_{ij} = \\
{\frac{{ \sum\limits_{s = 1}^{p}\sum\limits_{t =
1}^{n} \sum\limits_{l = 1}^{q} {\sum\limits_{\beta \in
J_{s_{1},\,m} {\left\{ {i} \right\}}} {{\rm{cdet}} _{i} \left(
{\left( {{\rm {\bf A}}{\bf W}_{1}} \right)^{k_{1}+2}_{\,. \,i}
\left( {\rm {\bf e}}_{.\,t} \right)} \right) {\kern
1pt} _{\beta} ^{\beta} } } }\hat{
d_{ts}}{\bar u}_{sl}{\sum\limits_{\alpha \in I_{s_{2},q} {\left\{ {j}
\right\}}} {{\rm{rdet}} _{j} \left( {({\bf W}_{2}{\rm {\bf B}} )^{ k_{2}+2}_{j\,.\,} ({\rm {\bf e}}_{l.} )} \right)\,_{\alpha} ^{\alpha}
} } }{{{{\sum\limits_{\beta \in J_{s_{1},\,m}} {{\left|
{\left( { {\bf A}{\bf W}_{1}} \right)^{k_{1}+2}{\kern 1pt} _{\beta} ^{\beta}
}  \right|}}} }}{{{\sum\limits_{\alpha \in
I_{s_{2},\,q}}  {{\left| {\left({\bf W}_{2} {{\rm {\bf B}} } \right)^{k_{2}+2}{\kern
1pt}  _{\alpha} ^{\alpha} } \right|}}} }}}    }=\\
{\frac{{ \sum\limits_{t =
1}^{n} \sum\limits_{l = 1}^{q} {\sum\limits_{\beta \in
J_{s_{1},\,m} {\left\{ {i} \right\}}} {{\rm{cdet}} _{i} \left(
{\left( {{\rm {\bf A}}{\bf W}_{1}} \right)^{k_{1}+2}_{\,. \,i}
\left( {\rm {\bf e}}_{.\,t} \right)} \right) {\kern
1pt} _{\beta} ^{\beta} } } }
{\bar
d_{tl}}{\sum\limits_{\alpha \in I_{s_{2},q} {\left\{ {j}
\right\}}} {{\rm{rdet}} _{j} \left( {({\bf W}_{2}{\rm {\bf B}} )^{ k_{2}+2}_{j\,.\,} ({\rm {\bf e}}_{l.} )} \right)\,_{\alpha} ^{\alpha}
} } }{{{{\sum\limits_{\beta \in J_{s_{1},\,m}} {{\left|
{\left( { {\bf A}{\bf W}_{1}} \right)^{k_{1}+2}{\kern 1pt} _{\beta} ^{\beta}
}  \right|}}} }}{{{\sum\limits_{\alpha \in
I_{s_{2},\,q}}  {{\left| {\left({\bf W}_{2} {{\rm {\bf B}} } \right)^{k_{2}+2}{\kern
1pt}  _{\alpha} ^{\alpha} } \right|}}} }}}    }.
\end{multline}
Denote by
\begin{multline*}
 d^{{\rm {\bf A}}}_{il}:= \\
{\sum\limits_{\beta \in
J_{s_{1},\,m} {\left\{ {i} \right\}}} {{\rm{cdet}} _{i} \left(
{\left( {{\rm {\bf A}}{\bf W}_{1}} \right)^{k_{1}+2}_{\,. \,i} \left( {{\bf {\bar d}}}_{.\,l}
\right)} \right){\kern 1pt} _{\beta} ^{\beta} } }= \sum\limits_{t
= 1}^{n} {\sum\limits_{\beta \in
J_{s_{1},\,m} {\left\{ {i} \right\}}} {{\rm{cdet}} _{i} \left(
{\left( {{\rm {\bf A}}{\bf W}_{1}} \right)^{k_{1}+2}_{\,. \,i} \left( {\rm {\bf e}}_{.\,t}
\right)} \right){\kern 1pt} _{\beta} ^{\beta} } }
\bar{d}_{tl}
\end{multline*}
the $l$-th component  of a row-vector ${\rm {\bf d}}^{{\rm {\bf
A}}}_{i\,.}= (d^{{\rm {\bf A}}}_{i1},...,d^{{\rm {\bf A}}}_{iq})$
for all $l=\overline{1,q}$. Substituting it in (\ref{eq:x_ij}),
we have
\[x_{ij} ={\frac{{ \sum\limits_{l = 1}^{q}
 d^{{\rm {\bf A}}}_{il}
}{\sum\limits_{\alpha \in I_{s_{2},q} {\left\{ {j}
\right\}}} {{\rm{rdet}} _{j} \left( {({\bf W}_{2}{\rm {\bf B}} )^{ k_{2}+2}_{j\,.\,} ({\rm {\bf e}}_{l.} )} \right)\,_{\alpha} ^{\alpha}
} }
}{{{{\sum\limits_{\beta \in J_{s_{1},\,m}} {{\left|
{\left( { {\bf A}{\bf W}_{1}} \right)^{k_{1}+2}{\kern 1pt} _{\beta} ^{\beta}
}  \right|}}} }}{{\sum\limits_{\alpha \in
I_{s_{2},\,q}}  {{\left| {\left({\bf W}_{2} {{\rm {\bf B}} } \right)^{k_{2}+2}{\kern
1pt}  _{\alpha} ^{\alpha} } \right|}}} }}}.
\]
Since $\sum\limits_{l = 1}^{q}
 d^{{\rm {\bf A}}}_{il}{\rm {\bf e}}_{l.}={\rm {\bf
d}}^{{\rm {\bf A}}}_{i\,.}$, then it follows (\ref{eq:d^A}).

If we denote by
\begin{multline}\label{eq:d^B_den}
 d^{{\rm {\bf B}}}_{tj}:= \\
\sum\limits_{l= 1}^{q}{\bar
d_{tl}}{\sum\limits_{\alpha \in I_{s_{2},q} {\left\{ {j}
\right\}}} {{\rm{rdet}} _{j} \left( {({\bf W}_{2}{\rm {\bf B}} )^{ k_{2}+2}_{j\,.\,} ({\rm {\bf e}}_{l.} )} \right)_{\alpha} ^{\alpha}
} }={\sum\limits_{\alpha \in I_{s_{2},q} {\left\{ {j}
\right\}}} {{\rm{rdet}} _{j} \left( {({\bf W}_{2}{\rm {\bf B}} )^{ k_{2}+2}_{j\,.\,} ( {\bar {\bf
d}}_{t.} )} \right)_{\alpha} ^{\alpha}
} }
\end{multline}

\noindent the $t$-th component  of a column-vector ${\rm {\bf
d}}^{{\rm {\bf B}}}_{.\,j}= (d^{{\rm {\bf B}}}_{1j},...,d^{{\rm
{\bf B}}}_{nj})^{T}$ for all $t=\overline{1,n}$ and substituting it
in (\ref{eq:x_ij}), we obtain
\[x_{ij} ={\frac{{  \sum\limits_{t= 1}^{n}
{\sum\limits_{\beta \in
J_{s_{1},\,m} {\left\{ {i} \right\}}} {{\rm{cdet}} _{i} \left(
{\left( {{\rm {\bf A}}{\bf W}_{1}} \right)^{k_{1}+2}_{\,. \,i} \left( {\rm {\bf e}}_{.\,t}
\right)} \right){\kern 1pt} _{\beta} ^{\beta} } } }\,\,d^{{\rm {\bf
B}}}_{tj} }{{{\sum\limits_{\beta \in J_{r_{1},\,n}} {{\left|
{\left( {{\rm {\bf A}}^{ *} {\rm {\bf A}}} \right){\kern 1pt}
_{\beta} ^{\beta} }  \right|}}} }{{\sum\limits_{\alpha \in
I_{r_{2},p}} {{\left| {\left( {{\rm {\bf B}}{\rm {\bf B}}^{ *} }
\right) _{\alpha} ^{\alpha} } \right|}}} }}    }.
\]
Since $\sum\limits_{t = 1}^{n}{\rm {\bf e}}_{.t}
 d^{{\rm {\bf B}}}_{tj}={\rm {\bf
d}}^{{\rm {\bf B}}}_{.\,j}$, then it follows (\ref{eq:d^B}).
$\Box$

\section{Examples}
In this section, we give examples to illustrate our results.

1. Let
us consider the matrix equation

\begin{equation}\label{ex:ax}
{\bf W}{\bf A}{\bf W}{\bf X}={\bf D}\\
\end{equation}
with the restricted conditions (\ref{eq:ax_restr}), where
\[{\bf A}=\begin{pmatrix}
  0 & i & 0 \\
  k & 1 & i \\
 1 & 0 & 0\\
  1 & -k & -j
\end{pmatrix},\,\, {\bf W}=\begin{pmatrix}
 k & 0 & i & 0 \\
 -j & k & 0 & 1 \\
 0 & 1 & 0 & -k
\end{pmatrix}, {\bf D}=\begin{pmatrix}
 k &   i \\
 i & - j  \\
 1 & -i
\end{pmatrix}.\]
Then
\[{\bf V}={\bf A}{\bf W}=\begin{pmatrix}
 -k & -j & 0 & i \\
 -1-j & i+k & j & 1+j \\
 k & 0 & i & 0\\
 -i+k & 1-j & i & i-k
\end{pmatrix},\,\,{\bf U}={\bf W}{\bf A}=\begin{pmatrix}
  i & j & 0 \\
  0 & k & 0 \\
 0 & 0 & 0
\end{pmatrix},\]
and $\rank{\bf W}=3$, $\rank{\bf V}=3$, $\rank{\bf V}^{3}=\rank {\bf V}^{2}=2$, $\rank{\bf U}^{2}=\rank {\bf U}=2$.
So, $ Ind\,{\bf V}=2$, $ Ind\,{\bf U}=1$, and $k= {\rm max}\{Ind({\bf A}{\bf W}), Ind({\bf W}{\bf A})\}=2$.

 We shall find the  W-weighted Drazin inverse solution of (\ref{ex:ax}) by it's determinantal representation (\ref{eq:cr_rul_ax1}).
We have \[{\bf U}^{2}=\begin{pmatrix}
  -1 & i+k & 0 \\
  0 & -1 &  \\
 0 & 0 & 0
\end{pmatrix},\,\,\,{\bf U}^{5}=\begin{pmatrix}
  i & 2+3j & 0 \\
  0 & k &  \\
 0 & 0 & 0
\end{pmatrix},\]\[({\bf U}^{5})^{*}=\begin{pmatrix}
  -i & 0 & 0 \\
  2-3j & -k &  \\
 0 & 0 & 0
\end{pmatrix},\left({\bf U}^{5}\right)^{*}{\bf U}^{5}=\begin{pmatrix}
  1 & -2i-3k & 0 \\
  2i+3k & 14 &  \\
 0 & 0 & 0
\end{pmatrix},\] \[
\hat{{\rm {\bf D}}}=({\bf U}^{ 5})^{*}{\bf U}^{2}{\bf D}=\begin{pmatrix}
  i-j-k & -j \\
  1+3i+6j-2k & 4i-2k  \\
 0 & 0
\end{pmatrix},{\bf W}^{*}=\begin{pmatrix}
 -k & j & 0  \\
 0 & -k & 1 \\
 -i & 0 &  0\\
 0 & 1 & k
\end{pmatrix},\]\[
{\bf W}^{*}{\bf W}=\begin{pmatrix}
 2 & i & -j & j  \\
 -i & 2 & 0 & -2k \\
 j & 0 &  1 & 0\\
  -j & 2k & 0 & 2 \\
\end{pmatrix},
\hat{{\rm {\bf W}}}={\bf W}^{*}{\bf U}^{2}=\begin{pmatrix}
 -k & 1-2j & 0  \\
 0 & i+k & 0 \\
 i & 1+j &  0\\
 0 & -1 & 0
\end{pmatrix}.
\]
Since by (\ref{eq:cr_rul_ax1})
\begin{multline*}
 {x}_{11}=\\
\frac{\sum\limits_{t = 1}^{3}
{{\sum\limits_{\beta \in I_{3,\,4} {\left\{ {1}
\right\}}} {{\rm{cdet}} _{1} \left( {\left(   {\bf W}^{*} {\bf W}\right)_{.1 } (\hat{{\rm {\bf w}}}_{.t})} \right) {\kern 1pt} _{\beta} ^{\beta} } }
}
{{\sum\limits_{\beta \in J_{2,\,3} {\left\{ {t}
\right\}}} {{\rm{cdet}} _{t} \left( {\left( { \left({\bf U}^{ 5} \right)^{*}{\bf U}^{ 5}
} \right)_{.t } (\hat{{\rm {\bf d}}}_{.1})} \right) {\kern 1pt} _{\beta} ^{\beta} } }
}}
{{\sum\limits_{\beta \in J_{3,\,4}} {{\left| {\left(  {\bf W}^{*}{\bf W}\right){\kern 1pt} _{\beta} ^{\beta}
}  \right|}}}{\sum\limits_{\beta \in J_{2,\,
3}} {{\left| {\left( { \left({\bf U}^{ 5} \right)^{*}{\bf U}^{ 5}
} \right){\kern 1pt} _{\beta} ^{\beta}
}  \right|}}}},\end{multline*}
where
\begin{multline*}{{\sum\limits_{\beta \in I_{3,\,4} {\left\{ {1}
\right\}}} {{\rm{cdet}} _{1} \left( {\left(   {\bf W}^{*} {\bf W}\right)_{.1 } (\hat{{\rm {\bf w}}}_{.1})} \right) {\kern 1pt} _{\beta} ^{\beta} } }
}=\\{\rm{cdet}} _{1}\begin{pmatrix}
 k & i & -j  \\
 0 & 2 & 0 \\
 i & 0 &  1
\end{pmatrix}+{\rm{cdet}} _{1}\begin{pmatrix}
 k & i & j  \\
 0 & 2 & -2k \\
 0 & 2k &  1
\end{pmatrix}+{\rm{cdet}} _{1}\begin{pmatrix}
 k & -j & j  \\
 i & 1 & 0 \\
 0 & 0 &  2
\end{pmatrix}=0, \end{multline*}
\begin{multline*}{{\sum\limits_{\beta \in I_{3,\,4} {\left\{ {1}
\right\}}} {{\rm{cdet}} _{1} \left( {\left(   {\bf W}^{*} {\bf W}\right)_{.1 } (\hat{{\rm {\bf w}}}_{.2})} \right) {\kern 1pt} _{\beta} ^{\beta} } }
}=-2j,\\
{{\sum\limits_{\beta \in I_{3,\,4} {\left\{ {1}
\right\}}} {{\rm{cdet}} _{1} \left( {\left(   {\bf W}^{*} {\bf W}\right)_{.1 } (\hat{{\rm {\bf w}}}_{.3})} \right) {\kern 1pt} _{\beta} ^{\beta} } }}=0,
{\sum\limits_{\beta \in J_{3,\,4}} {{\left| {\left(  {\bf W}^{*}{\bf W}\right){\kern 1pt} _{\beta} ^{\beta}
}  \right|}}}=2,\end{multline*}
and
\begin{multline*}\sum\limits_{\beta \in J_{2,\,3} {\left\{ {1}
\right\}}}{{\rm{cdet}} _{1} \left( {\left( { \left({\bf U}^{ 5} \right)^{*}{\bf U}^{ 5}
} \right)_{.1 } (\hat{{\rm {\bf d}}}_{.1})} \right) {\kern 1pt} _{\beta} ^{\beta} }=\\{\rm{cdet}} _{1}\begin{pmatrix}
 i-j-k & -2i-3k   \\
 1+3i+6j-2k & 14
\end{pmatrix}+{\rm{cdet}} _{1}\begin{pmatrix}
 i-j-k & 0   \\
 0 & 0
\end{pmatrix}=-2i-j-k,\end{multline*}
\begin{multline*}\sum\limits_{\beta \in J_{2,\,3} {\left\{ {2}
\right\}}}{{\rm{cdet}} _{2} \left( {\left( { \left({\bf U}^{ 5} \right)^{*}{\bf U}^{ 5}
} \right)_{.2 } (\hat{{\rm {\bf d}}}_{.1})} \right) {\kern 1pt} _{\beta} ^{\beta} }=j,\\\sum\limits_{\beta \in J_{2,\,3} {\left\{ {3}
\right\}}}{{\rm{cdet}} _{3} \left( {\left( { \left({\bf U}^{ 5} \right)^{*}{\bf U}^{ 5}
} \right)_{.3 } (\hat{{\rm {\bf d}}}_{.1})} \right) {\kern 1pt} _{\beta} ^{\beta} }=0, {\sum\limits_{\beta \in J_{2,\,
3}} {{\left| {\left( { \left({\bf U}^{ 5} \right)^{*}{\bf U}^{ 5}
} \right){\kern 1pt} _{\beta} ^{\beta}
}  \right|}}}=1,\end{multline*}
then
\begin{align*}
 {x}_{11}=&
\frac{0\cdot (-2i-j-k)+(-2j)\cdot j+0\cdot 0}{2\cdot 1 }=1,\\
 {x}_{12}=&\frac{0\cdot (-2+2j)+(-2j)\cdot i+0\cdot 0}{2\cdot 1 }=k,\\
 {x}_{21}=&\frac{2j\cdot (-2i-j-k)+(10i-4k)\cdot j+0\cdot 0}{2\cdot 1 }=1+i+7k,\\
 {x}_{22}=&\frac{2j\cdot (-2+2j)+(10i-4k)\cdot i+0\cdot 0}{2\cdot 1 }=-7-4j,\\
 {x}_{31}=&\frac{10i\cdot (-2i-j-k)+j\cdot j+0\cdot 0}{2\cdot 1 }=9.5+5j-5k,\\
 {x}_{32}=&\frac{10i\cdot (-2+2j)+j\cdot i+0\cdot 0}{2\cdot 1 }=-10i+9.5k,
\end{align*}
We finally get,
 \begin{equation}
\label{ex:wdr_rep1}{\rm {\bf X}}= \begin{pmatrix}
  1 &k \\
  1+i+7k &-7-4j   \\
9.5+5j-5k &-10i+9.5k
\end{pmatrix}.\end{equation}

2. Let now we consider the matrix equation
\begin{equation}\label{ex:axb}
{\bf W}_{1}{\bf A}{\bf W}_{1}{\bf X}{\bf W}_{2}{\bf B}{\bf W}_{2}={\bf D},
\end{equation}
 with the constraints (\ref{eq:axbr}), where

\[{\bf A}=\begin{pmatrix}
 k & 0 & i & 0 \\
 -j & k & 0 & 1 \\
 0 & 1 & 0 & -k
\end{pmatrix},\,{\bf W}_{1}=\begin{pmatrix}
 k &- j & 0\\
 0 &k&1\\
 i & 0 & 0 \\
 0 & 1 & -k
\end{pmatrix}, {\bf W}_{2}=\begin{pmatrix}
k &-i \\
 j &0\\
 0 & 1
\end{pmatrix},\] \[{\bf B} =\begin{pmatrix}
 k & j & 0\\
 j &0&1
\end{pmatrix}, {\bf D}=\begin{pmatrix}
 i &-1 \\
 k &0\\
 0 &j\\
 -1& 0
\end{pmatrix}.\]
Since the following matrices are Hermitian
\[{\bf V}={\bf A}{\bf W}_{1}=\begin{pmatrix}
 -2 & i & 0  \\
 -i & -1 & 0\\
 0 & 0 &  0
\end{pmatrix},\,\,{\bf U}={\bf W}_{2}{\bf B}=\begin{pmatrix}
 0 & -i&-i   \\
 i&-1  & 0\\
 i&0  & -1
\end{pmatrix},\]
then we can find the  W-weighted Drazin inverse solution of (\ref{ex:axb}) by it's determinantal representation (\ref{eq:d^B}).

We have
\begin{gather*}k_{1}=\max\,\{Ind({\bf A}{\bf W}_{1}), Ind\,({\bf W}_{1}{\bf A})\}=1,\\
k_{2}=\max\,\{Ind({\bf B}{\bf W}_{2}),Ind\,({\bf W}_{2}{\bf B})\}=1,\end{gather*}
and $s_{1}=\rank\,( {\rm {\bf A}}{\bf W}_{1}) =2$, $s_{2}=\rank\,({\bf W}_{2} {\rm {\bf B}}) =2$. Since
\[( {\rm {\bf A}}{\bf W}_{1})^{3}=
\begin{pmatrix}
 -13 &8i & 0 \\
 -8i & -5 & 0\\
 0 & 0 &  0
\end{pmatrix},
({\bf W}_{2} {\rm {\bf B}})^{3}=
\begin{pmatrix}
0 &-3i & -3i \\
 3i & -3 & 0\\
 3i & 0 &  3
 \end{pmatrix},\]
then \[{\sum\limits_{\beta \in J_{2,\,3}} {{\left|
{\left( { {\bf A}{\bf W}_{1}} \right)^{3}{\kern 1pt} _{\beta} ^{\beta}
}  \right|}}}=1, \sum\limits_{\alpha \in I_{2,\,3}}
  {{\left| {\left({\bf W}_{2} {{\rm {\bf B}} } \right)^{3}}\, _{\alpha} ^{\alpha}  \right|}}=-27.\]
We have
\[    {\bf \bar{D}}={\bf A}{\bf W}_{1}{\bf A}{\bf D}{\bf B}{\bf W}_{2}{\bf B}=\begin{pmatrix}
 2i+j &-7+k & -5+2k \\
   -1+k & -5i-j & -4i-2j\\
   0 & 0 &  0
   \end{pmatrix},\]

 By (\ref{eq:def_d^B_m}), we can get \[{\bf d}_{.1}^{
 {\bf B}}=\begin{pmatrix}
   36i-9j \\
   -27-9k \\
   0
 \end{pmatrix},\,\,\,\,\,{\bf d}_{.2}^{{\bf  B}}=\begin{pmatrix}
   -27 \\ -18i \\
   0
 \end{pmatrix},\,\,\,\,\,{\bf d}_{.3}^{{\bf  B}}=\begin{pmatrix}
   9-9k \\ 9i+3j \\
   0
 \end{pmatrix}.\]
Since
\[{\left( {{\rm {\bf A}}{\bf W}_{1}} \right)^{3}_{\,. \,1} \left( {{{\rm {\bf d}}}_{.\,1}^{{\rm {\bf B}}}}
\right)}=
\begin{pmatrix}
  36i-9j & 8i & 0 \\
  -27-9k & -5 & 0 \\
 0 & 0 & 4
\end{pmatrix},\] then finally we obtain
\[
  x_{11} = {\frac{{\sum\limits_{\beta \in
J_{2,\,3} {\left\{ {1} \right\}}} {{\rm{cdet}} _{1} \left(
{\left( {{\rm {\bf A}}{\bf W}_{1}} \right)^{3}_{\,. \,1} \left( {{{\rm {\bf d}}}_{.\,1}^{{\rm {\bf B}}}}
\right)} \right){\kern 1pt} _{\beta} ^{\beta} } }}{{{{\sum\limits_{\beta \in J_{2,\,3}} {{\left|
{\left( { {\bf A}{\bf W}_{1}} \right)^{3}{\kern 1pt} _{\beta} ^{\beta}
}  \right|}}} }}{{\sum\limits_{\alpha \in
I_{2,\,3}}  {{\left| {\left({\bf W}_{2} {{\rm {\bf B}} } \right)^{3}{\kern
1pt}  _{\alpha} ^{\alpha} } \right|}}} }}}=\frac{36i-27j}{-27}=\frac{-4i+3j}{3},\]
Similarly,
\begin{gather*}
x_{12} =\frac{{\rm{cdet}} _{1}\begin{pmatrix}
  -27 & 8i  \\
  -18i & -5
\end{pmatrix}}{-27}=\frac{1}{3},\,\,
x_{13} =\frac{{\rm{cdet}} _{1}\begin{pmatrix}
  9-9k & 8i  \\
  9i-3j & -5
\end{pmatrix}}{-27}=\frac{-9-7k}{9},\\
x_{21} =\frac{{\rm{cdet}} _{2}\begin{pmatrix}
  -13 & 36i-9j  \\
  -8i & -27-9k
\end{pmatrix}}{-27}=\frac{-7-5k}{3},\,\,x_{22} =\frac{{\rm{cdet}} _{2}\begin{pmatrix}
  -13 & -27  \\
  -8i & -18i
\end{pmatrix}}{-27}=\frac{-2i}{3},\\
x_{23} =\frac{{\rm{cdet}} _{2}\begin{pmatrix}
  -13 & -9-9k  \\
  -8i & 9i+3j
\end{pmatrix}}{-27}=\frac{15i-11j}{9},\,\,x_{31}=x_{32}=x_{33}=0.
\end{gather*}
So, the  W-weighted Drazin inverse solution of (\ref{ex:axb}) are
\[{\bf X}=\frac{1}{9}\begin{pmatrix}
 -12i+9j &3 & -9-7k \\
 -21-15k & -6i & 15i-11j\\
 0 & 0 &  0
\end{pmatrix}.\]
Note that we used Maple with the package CLIFFORD in the calculations.


\begin{thebibliography}{40}
\bibitem{ad}S. L. Adler, Quaternionic Quantum Mechanics and Quantum Fields. Oxford University
Press, New York, 1995.
\bibitem{gi} J.D. Gibbon, A quaternionic structure in the three-dimensional Euler and ideal magnetohydrodynamics equation. \emph{Physica D} 166 (2002) 17-28.
\bibitem{gi1}    J.D. Gibbon, D.D. Holm, R.M. Kerr, I. Roulstone, Quaternions and particle dynamics in the Euler fluid equations. \emph{Nonlinearity} 19 (2006) 1969--1983.

\bibitem{ha}A. Handson, H. Hui, Quaternion frame approach to streamline visualization. \emph{IEEE Tran.
Vis. Comput. Grap.} 1 (1995) 164-172.
\bibitem{st}B.L.Stevens,  F.L.Lewis, E.N.Johnson, Aircraft Control and Simulation: Dynamics, Controls Design, and Autonomous Systems. 3rd Edition, Wiley, 2015.

\bibitem{pe}A. Perez, J. M. McCarthy, Dual quaternion synthesis of constrained robotic systems. \emph{J. Mech. Des.} (2003) 126(3) 425--435.
\bibitem{is}T. Isokawa, T. Kusakabe, N. Matsui, F. Peper, Quaternion neural network and its application. \emph{Lecture Notes in Computer Science} 2774 (2003) 318-324.
\bibitem{pro} J. Pro\v{s}kov\v{a}, Description of protein secondary structure using dual quaternions. \emph{Journal
of Molecular Structure} 1076 (2014) 89-93.



\bibitem{wa1} S.F. Yuan, Q.W. Wang, X.F. Duan, On solutions of the quaternion matrix equation $AX=B$  and their applications in color image restoration. \emph{Appl. Math. Comput.} 221 (2013) 10--20.
\bibitem{wa2}Q. W. Wang, S. W. Yu, Extreme ranks of real matrices in solution of the
quaternion matrix equation $AXB = C$ with applications. \emph{Algebra Colloq.} 17(2)
(2010) 345--360.
\bibitem{yu}S. Yuan, A. Liao, Y. Lei, Least squares Hermitian solution of the matrix equation $(AXB,CXD)=(E,F)$ with the least norm over the skew field of quaternions. \emph{Math. Comput. Model.}  48 (2008) 91--100.
\bibitem{fe} L. G. Feng, W. Cheng, The solution set to the quaternion matrix equation
$AX -  \bar{X}B = 0$. \emph{Algebra Colloq.} 19(1) (2012) 175--180.
\bibitem{ji}T. S. Jiang, M. S. Wei, On a solution of the quaternion matrix equation
$X -  A \tilde{X}B = C$. \emph{Acta Math. Sin. (Engl. Ser.)} 21(3) (2005) 483-490.
\bibitem{son}C. Q. Song, G. L. Chen, X. D. Wang, On solutions of quaternion matrix
equations $XF -  AX = BY$ and $XF  - A \tilde{X} = BY$.  \emph{Acta Math. Sci.} 32(5) (2012) 1967--1982.
 \bibitem{yu1}S.Yuan,  Q. W. Wang,  Two special kinds of least squares solutions for the quaternion matrix equation
$AXB+CXD=E$. \emph{Electron. J. Linear Algebra} 23 (2012) 57--74.
\bibitem{zha} F. Zhang, M. Wei, Y. Lia, J. Zhao, Special least squares solutions of the quaternion matrix equation AX=B with applications. \emph{Appl. Math. Comput.} 270 (2015) 425--433.



\bibitem{cl}
R.E. Cline, T.N.E. Greville, A Drazin inverse for rectangular matrices. \emph{Linear Algebra
Appl.} 29 (1980) 53--62.
 \bibitem{ca} S. L. Campbell, C. D. Meyer Jr., N. J. Rose, Applications of the Drazin inverse to linear systems of differential equations with singular constant coefficients \emph{SIAM J. Appl. Math}. 31 (1976)  411-425.
\bibitem{camp}S.L. Campbell, Comments on 2-D descriptor systems. \emph{Automatica} 27(1) (1991) 189--192.
\bibitem{ca1} S. L. Campbell, C.D. Meyer, Generalized inverse of linear transformations. Corrected reprint of the 1979 original. Dover Publications, Inc., New York, 1991.
\bibitem{sp}D. J. Spitzner, T. R Boucher, Asymptotic variance of functionals of discrete-time Markov chains via the Drazin inverse. \emph{Elect. Comm. in Probab.} 12 (2007) 120--133.
\bibitem{ka}T. Kaczorek, Application of the Drazin inverse to the analysis of
descriptor fractional discrete–time linear
systems with regular pencils. \emph{Int. J. Appl. Math. Comput. Sci.} 23(1) (2013) 29--33.
\bibitem{ze} Z. Al-Zhour, A. Kili\c{c}man, M. H. Abu Hassa, New representations for weighted
Drazin inverse of matrices. \emph{Int. Journal of Math. Analysis} 1(15) (2007) 697--708.
\bibitem{nik} M. Nikuie, M. Z. Ahmad, New Results on the W-weighted Drazin Inverse. \emph{AIP Conference Proceedings} 1602 (2014) 157.
\bibitem{nik1}M. Nikuie, Singular fuzzy linear systems. \emph{App. Math. and Comp. Intel.} 2(2) (2013) 157-168.
 \bibitem{wang}   G. Wang, J. Sun. A Cramer rule for solution of the general restricted matrix equation. \emph{Appl. Math. Comput.} 154 (2004) 415–422.

\bibitem{wei}
Y. Wei. A characterization for the W-weighted Drazin inverse and a Cramer rule for the
W-weighted Drazin inverse solution. \emph{Appl. Math. Comput.} 125 (2002) 303--310.







\bibitem{song_ela}
G. Song, Characterization of the W-weighted Drazin inverse over the quaternion
skew field with applications. \emph{Electron. J. Linear Algebra} 26 (2013) 1--14.
   \bibitem{ky11} I. Kyrchei, Determinantal representations of the W-weighted Drazin inverse over the quaternion skew field. \emph{Appl. Math. Comput.} 264 (2015) 453--465.


 \bibitem{ky1} I. Kyrchei, Cramer's rule for quaternion
 systems of linear equations. \emph{J. Math. Sci.} 155(6) (2008) 839--858.

 \bibitem{ky2}    I. Kyrchei, The theory of the column and row determinants in a quaternion linear algebra.  In: Albert R. Baswell (Eds.), Advances in Mathematics Research 15,  Nova Sci. Publ., New York, pp. 301--359, 2012.
   \bibitem{ky3}  I. Kyrchei, Cramer's rule for some quaternion matrix equations. \emph{Appl. Math. Comput.} 217(5) (2010) 2024--2030.





     \bibitem{ky4} I. Kyrchei, Determinantal representations of the Moore-Penrose inverse over the quaternion skew field and corresponding Cramer's rules. \emph{Linear Multilinear A.} 59 (2011) 413--431.
\bibitem{ky7} I. Kyrchei, Determinantal representation of the Moore–Penrose inverse matrix over the quaternion skew field. \emph{J. Math. Sci.} 180(1) (2012) 23--33.
   \bibitem{ky6}    I. Kyrchei, Explicit representation formulas for the minimum norm least squares solutions of some quaternion matrix equations.  \emph{Linear Algebra Appl.}  438 (2013) 136--152.

   \bibitem{ky5} I. Kyrchei, Determinantal representations of the Drazin inverse over the quaternion skew field with applications to some matrix equations. \emph{Appl. Math. Comput.} 238 (2014) 193--207.

\bibitem{song3}G. Song,   Determinantal representation of the generalized inverses over the quaternion skew field with applications. \emph{Appl. Math. Comput.} 219 (2012) 656--667.
\bibitem{song5}G. Song, Bott-Duffin inverse over the quaternion skew field with applications. \emph{J. Appl. Math. Comput.} 41 (2013) 377--392.
\bibitem{song1}G. Song, Q. Wang, H. Chang, Cramer rule for the unique solution of restricted matrix
equations over the quaternion skew field. \emph{Comput. Math. Appl.} 61 (2011) 1576--1589.




 \bibitem{song6}G. Song, X. Wang, X. Zhang, On solutions of the generalized Stein quaternion matrix equation. \emph{J. Appl. Math. Comput.} 43 (2013) 115--131.

     \bibitem{ra} V. Rako\v{c}evi\'{c}, Y. Wei, A weighted Drazin inverse and applications. \emph{Linear Algebra
Appl.} 350 (2002) 25--39.

\bibitem{wei1} Y. Wei. Integral representation of the W-weighted Drazin inverse. \emph{Appl. Math. Comput.} 144 (2003) 3--10.





\end{thebibliography}
\end{document}